\documentclass[11pt]{amsart}

\usepackage{amssymb, amsmath, eufrak, graphics}

\numberwithin{equation}{section}

\newcommand{\sgn}{sgn}
\newcommand{\Arg}{\text{Arg}}

\newcommand{\TM}{T_1 \mathcal{M}}

\newcommand{\KKK}{\text{\rm K}}
\newcommand{\LLL}{\Lambda}
\newcommand{\Ll}{L}

\newcommand{\PSLqspace}{\Gamma  \setminus \text{PSL}(2,\R)}
\newcommand{\Tqspace}{\Gamma \setminus T_1 \mathcal{H}}

\newcommand{\minspecX}{X}
\newcommand{\minexp}{\nu}
\newcommand{\minss}{\mathcal{S}_{\Gamma}}
\newcommand{\ii}{{i}}
\newcommand{\emdash}{\: - \negthickspace \negthickspace - \:}

\newcommand{\sfrac}[2]{{\textstyle \frac {#1}{#2}}}
\newcommand{\frh}{\textstyle \frac 12}
\newcommand{\DD}{\mathcal{D}}
\newcommand{\matr}[4]{\left( \begin{smallmatrix} #1 & #2 \\ #3 & #4
\end{smallmatrix} \right) }

\newcommand{\FF}{\mathcal{F}}

\newcommand{\ygg}{\mathcal{Y}_{\Gamma}}

\newcommand{\CC}{\mathcal{C}}
\newcommand{\F}{\mathcal{F}}

\newcommand{\PSL}{\text{PSL}}
\newcommand{\SL}{\text{SL}}

\newcommand{\qspace}{\Gamma \setminus \mathcal{H}}

\newcommand{\Z}{{\mathbb Z}}

\newcommand{\R}{{\mathbb R}}
\newcommand{\C}{{\mathbb C}}
\newcommand{\M}{\mathcal{M}}
\newcommand{\HH}{\mathcal{H}}

\newcommand{\tim}{\text{ Im }}
\newcommand{\tre} {\text{Re }}

\newtheorem{theorem}{$\quad$Theorem}
\newtheorem{proposition}{$\quad$Proposition}[section]
\newtheorem{corollary}[proposition]{$\quad$Corollary}
\newtheorem{lemma}[proposition]{$\quad$Lemma}

\theoremstyle{definition}
\newtheorem{remark}[proposition]{$\quad$Remark}

\newenvironment{Proof}{{\bf Proof. }}{$\square$\newline}

\newcounter{bibno}

\begin{document}


\title[Equidistribution along closed horocycles]{Equidistribution 
of Kronecker sequences along closed horocycles}

\author{Jens Marklof and Andreas Str\"ombergsson}



\thanks{J.M. has been supported by an EPSRC Advanced
Research Fellowship, the Nuffield Foundation (Grant NAL/00351/G),
the Royal Society (Grant 22355), and
the European Commission 
under the Research Training Network (Mathematical Aspects of Quantum Chaos) 
HPRN-CT-2000-00103 of the IHP Programme.}

\thanks{A.S. has been fully supported by 
the European Commission 
under the Research Training Network (Mathematical Aspects of Quantum Chaos) 
HPRN-CT-2000-00103 of the IHP Programme.}

\begin{abstract}
It is well known that (i) for every irrational number $\alpha$
the Kronecker sequence $m\alpha$ ($m=1,\ldots,M$)
is equidistributed modulo one in the limit $M\to\infty$, and (ii)
closed horocycles of length $\ell$ become equidistributed
in the unit tangent bundle $\TM$ of a hyperbolic surface $\M$
of finite area, as $\ell\to\infty$.   
In the present paper 
both equidistribution problems are studied simultaneously:
we prove that for any constant ${\minexp}>0$
the Kronecker sequence embedded in $\TM$ along a long closed horocycle
becomes equidistributed in $\TM$ for almost all $\alpha$, provided
that $\ell=M^{\minexp}\to\infty$. 
This equidistribution result holds in fact under explicit
diophantine conditions on $\alpha$ (e.g., for $\alpha=\sqrt 2$)
provided that $\nu<1$, or $\nu<2$ with additional assumptions on
the Fourier coefficients of certain automorphic forms. Finally, we 
show that for $\nu=2$, our equidistribution theorem implies
a recent result of Rudnick and Sarnak on the uniformity of the
pair correlation density of the sequence $n^2 \alpha$ modulo one.
\end{abstract}

\maketitle

\section{Introduction}

New developments in the ergodic theory of unipotent flows have,
in the past, led to
the solution of important problems in number theory. 
A famous example is Margulis'
proof of the Oppenheim conjecture on values of quadratic forms,
and the subsequent proof of a quantitative version \cite{Eskin98}.
The latter crucially uses 
Ratner's Theorem \cite{Ratner91,Ratner91b}, 
which provides a complete description of 
all invariant ergodic measures of a unipotent flow. 

The present work is motivated by recent studies of the local spacing
distributions of the sequence $n^2\alpha$ modulo one ($n=1,2,3,\ldots$),
which are conjectured to coincide---for generic values of $\alpha$---with 
the spacing distribution
of independent random variables from a Poisson process 
\cite{RS,RSZ}. We will show that the convergence of the pair correlation
density of $n^2\alpha$ mod one
to the Poisson answer is implied by the equidistribution of Kronecker
sequences along closed unipotent orbits (horocycles)
in the unit tangent bundle $\TM$
of a non-compact hyperbolic surface $\M$ with finite area.
The major difficulty in proving equidistribution 
for the unipotent cascades considered here is that it is unknown if 
all possible limit measures are necessarily invariant under a unipotent
element, and hence Ratner's theory cannot be applied. 

To be more precise,
let us realize the hyperbolic surface $\M$ as a quotient
$\M=\qspace$, where $\Gamma$ is a cofinite Fuchsian group 
acting on the Poincar\'e upper half-plane 
$$
\HH=\{ x+\ii y \in \C \mid x\in\R, y>0\}
$$
with metric $ds = |dz|/y$, where $dz=dx+\ii dy$ is the complex line element.
We assume that $\M$ is non-compact, i.e., 
has at least one cusp. After a suitable coordinate transformation we may 
assume that one of the cusps lies at infinity, and that the corresponding
isotropy subgroup $\Gamma_{\infty}\subset\Gamma$ 
is generated by the translation $z \mapsto z+1$.

Let $T_1 \HH$ be the unit tangent bundle of $\HH$,
and denote its elements as $(z,\theta)$, with $z \in \HH$ and 
$\theta \in\R / 2\pi \Z$, 
where $\theta$ is an angular variable measured from the
vertical counterclockwise. The action of an element
$\gamma\in\Gamma$ on $T_1\HH$ is then given by
\begin{equation} \label{PSLACTIONDEF}
(z,\theta)\mapsto (\gamma z, \theta - 2\beta_\gamma(z))
\end{equation}
where
$$
\gamma=\begin{pmatrix} a & b \\ c & d \end{pmatrix},
\quad
\gamma z= \frac{a z +b}{c z +d},
\quad
\beta_\gamma(z)=\arg(cz+d).
$$
We have in particular $T_1 \M=\Tqspace$.

For any $y>0$, the curve $\{ (x + \ii y,0) \mid x \in \R \}$
is an example of an orbit of the horocycle flow on $T_1 \HH$.
By our assumption on $\Gamma_{\infty}$,
the above orbit will be closed in $T_1\M$, with length $y^{-1}$.
In the limit $y \to 0$, the orbit in fact becomes 
\textit{equidistributed} on $T_1 \M$ with respect to the 
Poincar\'e area $d\mu = dxdy/y^2$ 
times the uniform measure in the phase variable $\theta$:

\begin{theorem}[Sarnak \cite{Sarnak}] \label{sarnaks}
For any bounded continuous function $f: T_1\M \to \C$
we have
\begin{equation} 
\lim_{y\to 0}\int_0^1 f \left( x + {\ii y} ,0 \right) dx
= \langle f \rangle,
\end{equation}
were $\langle f \rangle$ denotes the average of $f$ over 
$T_1 \M$,
\begin{equation*}
\langle f \rangle = \frac 1{2\pi \mu(\F) } \int_{\M}
\int_0^{2\pi} f(z,\theta) \, d\theta \, d\mu(z).
\end{equation*}
\end{theorem}

The uniform measure on the closed orbit in Theorem \ref{sarnaks}
may in fact be replaced by uniform measures supported on 
comparatively small sub-arcs of the orbit \cite{Strhoro}.

Our investigation is concerned with the equidistribution of
the point set
\begin{equation} \label{POINTSET}
\Gamma \bigl\{ (m \alpha + {\ii y},0) \mid m=1,...,M \bigr\}
\end{equation}
as $M \to \infty$ and $y \to 0$. 
Clearly, it can be expected to be easier to 
establish equidistribution the faster the number of points 
$M$ grows relative to the hyperbolic length of the orbit, $y^{-1}$. 
In particular, when $M \gg y^{-1/{\minexp}}$
for some $0<{\minexp}<1$, we will see that equidistribution is a direct
consequence of the equidistribution of $m\alpha$ mod one and
Theorem \ref{sarnaks} (provided $\alpha$ is badly
approximable, cf.\ Remark \ref{SARNAKREMARK} below). 
The harder case is when $M$ is small compared with $y^{-1}$, and 
especially, as we shall see, when $M \ll y^{-1/2}$.

Our main result is the following.


\begin{theorem} \label{MAINTHM}
Fix ${\minexp}>0$.
Then there is a set $P=P(\Gamma,{\minexp}) \subset \R$ of full 
Lebesgue measure such
that for any $\alpha \in P$, any bounded continuous function
$f: T_1 \M \to \C$, and any constants $0<C_1<C_2$, we have
\begin{equation} \label{MAINTHMRES}
\frac 1M \sum_{m=1}^{M} f \left( m\alpha + {\ii y} ,0 \right) 
\to \langle f \rangle,
\end{equation}
uniformly as $M \to \infty$ and 
$C_1 M^{-{\minexp}} \leqq y \leqq C_2 M^{-{\minexp}}$.
\end{theorem}

In fact, we will prove a slightly stronger result than Theorem \ref{MAINTHM}
(cf.\ Theorem \ref{MAINTHM}' in Section \ref{MAINTHMPRSE}),
wherein the test function $f$ is allowed to be unbounded, satisfying a certain
growth condition in the cusps.
This will be important for the application to the distribution of $n^2 \alpha$ mod one. As we will show in Section \ref{APPLSECTION}, if relation
(\ref{MAINTHMRES}) holds for a specific $\alpha$, then 
the pair correlation density of $n^2\alpha$ mod one is uniform, i.e.,
coincides with the correlation density of independent random variables
from a Poisson process. Theorem \ref{MAINTHM}' therefore implies
the result by Rudnick and Sarnak \cite{RS} on $n^2 \alpha$.

\begin{remark}
For any given fixed $\delta \in \R$,
Theorem \ref{MAINTHM} remains true if we replace
$f(m\alpha + {\ii y},0)$ by $f(\delta+m\alpha + {\ii y},0)$ 
in \eqref{MAINTHMRES}.
\end{remark}

\begin{remark}
The conclusion in Theorem \ref{MAINTHM} is certainly
\textit{false} for all \textit{rational} $\alpha$.
For if $\alpha=p/q$ with $p,q \in \Z$, $q > 0$, then the set of points
$\Gamma \{(m\alpha+ \ii y,0) \mid m=1,...,M\}$ 
on $\Tqspace$ has cardinality 
$\leqq q$, for \textit{all} $M,y$.
In Section \ref{MOTEXSECT} we will give a much larger set of counter-examples
in the case where $\Gamma$ is a subgroup of $\PSL(2,\Z)$.
\end{remark}

We will also prove a theorem which gives an explicit diophantine condition on 
$\alpha$ ensuring \eqref{MAINTHMRES} to hold.
An irrational number $\alpha \in \R$ is said to be
of \textit{type} $K$ if there exists a constant $C>0$ such that
\begin{equation*}
\bigl | \alpha - \frac pq \bigr | > \frac C{q^K}
\end{equation*}
for all $p,q \in \Z$, $q>0$.
The smallest possible value of $K$ is $K=2$,
and it is well known that for any given $K>2$,
the set of $\alpha$'s of type $K$ is of full Lebesgue measure in $\R$.
Cf.,\ e.g.,\ \cite[\S 1]{Schmidt}.

We recall that under the above normalization of the cusp at $\infty$,
any Maass waveform $\phi(z)$ (of weight $0$) on $\qspace$  
has a Fourier expansion involving the Macdonald-Bessel function:
\begin{equation} \label{THM2PREL1}
\phi(z) = c_0 y^{\frac 12 - ir} +
\sum_{n \neq 0} c_n \sqrt y K_{\ii r}(2 \pi |n| y) e(nx),
\end{equation}
where $-\frac 14 - r^2 \leqq 0$ is the eigenvalue of
$\phi$.
Furthermore, if $\eta_1=\infty,\eta_2,...,\eta_{\kappa}$ are the
inequivalent cusps of $\Gamma$, then we have for
the Eisenstein series $E_k(z,s)$ ($k \in \{1,...,\kappa\}$)
associated to the cusp $\eta_k$:
\begin{align} \label{THM2PREL2}
E_k(z,s) = \delta_{k1} y^s + \varphi_{k1}(s) y^{1-s}
+ \sum_{n \neq 0} \psi_{n,k}(s) \sqrt y K_{s-\frac 12} (2\pi |n| y) 
e(nx).
\end{align}
We fix some $s_0 \in (\frac 12,1)$ such that $s=s_0$ is not a pole of
any of the Eisenstein series $E_k(z,s)$, $k = 1,...,\kappa$.
We will assume that 
$\beta \geqq 0$ is a number
such that the following holds, for each $\varepsilon>0$:

\begin{subequations}
\begin{align} \label{BETACOND:A} 
& {\textstyle \text{\parbox[t]{4.0in}{For each fixed Maass waveform as in
\eqref{THM2PREL1}, we have 
$c_n = O_{\varepsilon}(|n|^{\beta+\varepsilon})$ as $|n| \to \infty$.}}}
\\ \label{BETACOND:B} 
& {\textstyle \text{\parbox[t]{4.0in}{For each holomorphic cusp form $\phi(z) 
= \sum_{n=1}^{\infty} c_n e(nz)$ of even integer weight $m \geqq 2$
on $\qspace$,
we have $c_n = O_{\varepsilon}(|n|^{\frac{m-1}2 + \beta + \varepsilon})$
as $|n| \to \infty$.}}}
\\ \label{BETACOND:C} 
& {\textstyle \text{\parbox[t]{4.0in}{For each fixed 
$s \in \bigl( \frac 12 + \ii[0, \infty) \bigr)
\cup \{ s_0 \}$, and all
$k \in \{1,...,\kappa\}$, we have
$\psi_{n,k}(s) = O_{\varepsilon}(|n|^{\beta+\varepsilon})$ 
as $|n| \to \infty$.}}}
\end{align}
\end{subequations}

Thus, if $\Gamma$ is a congruence subgroup of $\text{PSL}(2,\Z)$, 
and we assume the Ramanujan conjecture for Maass waveforms 
to be true, then any $\beta \geqq s_0-\frac 12$ will work.
In particular, since $s_0$ can be taken arbitrarily close to $\frac 12$, 
we may then take $\beta>0$ arbitrarily small.
\textit{Unconditionally}, for $\Gamma$ a congruence subgroup,
we know that \eqref{BETACOND:A}-\eqref{BETACOND:C}
hold for $\beta \geqq \max(\frac 7{64},s_0-\frac 12)$,
by a recent result of Kim and Sarnak \cite{KS}.
(Recall here that the Ramanujan conjecture has been proved in the
holomorphic case \cite{Deligne}, 
i.e.\ \eqref{BETACOND:B} is known to be true with $\beta = 0$.
Also, on congruence subgroups, 
\eqref{BETACOND:C} can be shown to hold with $\beta = s_0-\frac 12$.)

For \textit{general} 
$\Gamma$ the situation is very different:
It follows from elementary bounds that
\eqref{BETACOND:A}-\eqref{BETACOND:C} hold for
$\beta=s_0$ (cf.\ Lemma \ref{EISENSTEINRS} below regarding
\eqref{BETACOND:C}).
Furthermore, according to Bernstein and Reznikov \cite{BR}, 
\eqref{BETACOND:A} holds for $\beta=\frac 13$ if we 
restrict ourselves to the case of Maass waveforms which are cusp forms.
Also, by Good \cite{Good}, \eqref{BETACOND:B} holds with $\beta=\frac 13$  
as long as the weight $m$ is $>2$.
It can be expected that \eqref{BETACOND:B} with $\beta=\frac 13$ should 
also be provable for $m=2$,
and that it should be possible to extend the methods in \cite{BR}
to the case of the Eisenstein series
and non-cuspidal Maass waveforms,
so as to prove \eqref{BETACOND:A} with $\beta=\frac 13$ and
\eqref{BETACOND:C} with some $\beta=\beta(s_0)>\frac 13$
such that $\beta(s_0) \to \frac 13$ if $s_0 \to \frac 12$.
(We are grateful to P. Sarnak for discussions on these matters.)

\begin{theorem} \label{MAINTHM2}
Let $s_0 \in (\frac 12,1)$ and $\beta > 0$ be as in 
\eqref{BETACOND:A}--\eqref{BETACOND:C}.
Let $\alpha \in \R$ be of type $K \geqq 2$,
and let ${\minexp}$ be a positive number satisfying
\begin{equation} \label{MAINTHM2ASS}
{\minexp} <
\begin{cases} \frac{2}{1+2\beta} 
& \text{if } \beta < \frac{3-K}{2(K-1)}
\\ \frac{2}{2K\beta+K-2} 
& \text{if } \frac{3-K}{2(K-1)} \leqq \beta < \frac 12
\\ \frac 2{2K+2\beta-3}
& \text{if } \frac 12 \leqq \beta .
\end{cases}
\end{equation}
Then for any bounded continuous function
$f: T_1\M \to \C$, and any constant $C_1>0$, we have
\begin{equation} \label{MAINTHM2CONCL}
\frac 1M \sum_{m=1}^{M} f \left( m\alpha + {\ii y} ,0 \right) 
\to \langle f \rangle,
\end{equation}
uniformly as $M \to \infty$, $y \to 0^+$ 
so long as $y \geqq C_1 M^{-{\minexp}}$.
\end{theorem}


Notice that the right hand side in \eqref{MAINTHM2ASS}
is a continuous function of $\beta$ and $K$, for
$\beta \geqq 0,$ $K \geqq 2$.
The restriction on $\nu$ in \eqref{MAINTHM2ASS}
is the best possible which can be obtained by our method
of using the absolute bounds in \eqref{BETACOND:A}-\eqref{BETACOND:C},
cf.\ Remark \ref{BESTPOSSIBLE} below.

\begin{corollary} \label{RPCOROLLARY}
Let $\Gamma$ be a congruence subgroup of $\PSL(2,\Z)$ and 
assume the Ramanujan conjecture for Maass waveforms on $\qspace$ to hold.
Then, if $\alpha$ is of type $K \geqq 2$,
\eqref{MAINTHM2CONCL} holds for any positive number
$\nu < \begin{cases} 2
& \text{if } K<3
\\ \frac{2}{K-2} 
& \text{if } K \geqq 3.
\end{cases}$
\end{corollary}

\begin{remark} 
In the case $K \geqq 3$, 
the bound $\nu < \frac 2{K-2}$ in Corollary \ref{RPCOROLLARY}
is actually the best possible restriction
on $\nu$, as follows from Proposition \ref{NEGCONGRPROP} below.
To be precise, let $\alpha$ be any irrational number
and let $K_0$ be the infimum of all
numbers $K$ such that $\alpha$ is of type $K$.
Then if $2 < K_0 < \infty$,
\eqref{MAINTHM2CONCL} is \textit{false}
for each $\nu > \frac 2{K_0-2}$.
Furthermore, if 
$\liminf_{q \to \infty} \bigl( \inf_{p \in \Z} q^{K_0} |\alpha-p/q|
\bigr)$ 
is finite and sufficiently small,
then \eqref{MAINTHM2CONCL} is also false for $\nu = \frac 2{K_0-2}$.
\end{remark}

\begin{remark} 
Unconditionally,
if $\Gamma$ is a congruence subgroup,
it follows from Theorem \ref{MAINTHM2} 
and \cite{KS} that if 
$2 \leqq K \leqq \frac{103}{39}$,
then \eqref{MAINTHM2CONCL} holds for any $\nu < \frac{64}{39}$.
\end{remark}

\begin{remark} \label{SARNAKREMARK}
As we have pointed out, for any group $\Gamma$, 
\eqref{BETACOND:A}--\eqref{BETACOND:C} hold
for $\beta = s_0$.
Hence, Theorem \ref{MAINTHM2}
implies that \eqref{MAINTHM2CONCL} holds whenever 
${\minexp} < (K-1)^{-1}$.
However, this fact can also be derived more straighforwardly
from Sarnak's Theorem \ref{sarnaks},
since we can prove directly that under the above conditions
the set of points  
$\{ (m \alpha + {\ii y},0) \mid m=1,...,M \}$ tends to become more and more
equidistributed \textit{along the closed horocycle}
$\{ (x + {\ii y},0) \mid x \in [0,1] \}$ (in the hyperbolic metric).
We give an outline of this argument at the end of Section \ref{PROOFOFREMARK}.
\end{remark}

\section*{Acknowledgment}

We would like to thank H. Furstenberg, G. Margulis, 
Z. Rudnick, \mbox{P. Sarnak} and D. Witte for 
useful and inspiring discussions.

\section{Spectral preliminaries}

In this section we recall some basic facts concerning the
spectral expansion of functions on $\Tqspace$,
and collect some useful bounds and formulas.
We start by introducing some notation which will be in force
throughout this paper.

We let $\Gamma\subset \PSL(2,\R)$ be a cofinite Fuchsian group such
that $\qspace$ has at least one cusp.

Concerning the cusps and the fundamental domain,
we will use the same notation as in
\cite[p.\ 268]{TR2}.
Specifically: we let $\F \subset \mathcal{H}$ be a 
canonical (closed) fundamental domain for $\Gamma \setminus
\mathcal{H}$, 
and let $\eta_1,...,\eta_{\kappa}$ 
(where $\kappa \geqq 1$) be the
vertices of $\F$
along $\partial \HH = \R \cup \{ \infty \}$.
Since $\F$ is canonical, $\eta_1,...,\eta_{\kappa}$ are
$\Gamma-$inequivalent.

For each $k \in \{ 1,...,\kappa \}$ we choose 
$N_k \in \text{PSL} (2,\R)$ such that 
$N_k(\eta_k)= \infty$ and such that the stabilizer
$\Gamma_{\eta_k}$ is $[ T_k ] ,$ where  $T_k := 
N_k^{-1} \matr 1{-1}01 N_k.$ 
Since $\FF$ is canonical, 
by modifying $N_k$ we can \textit{also} ensure that  
\begin{align} \label{CKRAV}
& N_k(\F) \: \bigcap \: \{ z \in \mathcal{H} \mid \tim z
\geqq B \}   
\\ \notag
& = 
\{ z \in \mathcal{H} \mid 0 \leqq \tre z \leqq 1 , \; \tim z 
\geqq B \}, 
\end{align}
for all $B$ large enough. 
We \textit{fix}, once and for all, a constant $B_0 > 1$ such that
\eqref{CKRAV} holds for all   
$B \geqq B_0$ and all $k \in \{1,...,\kappa\}.$ 


In line with what we stated in the 
introduction, we will make the assumption that
$\Gamma$ is normalized so that
\begin{equation} \label{NORMALIZEDCUSP}
\eta_1 = \infty, \qquad N_1 = \begin{pmatrix} 1 & 0 \\ 0 & 1 
\end{pmatrix}.
\end{equation}
This means that $\Gamma$ has one cusp located at $\infty$ of standard
width $z \mapsto z+1$.

We recall the definition of the 
\textit{invariant height function,} $\ygg(z)$:
\begin{equation} \label{YGGDEF}
\ygg (z) 
= \sup \Bigl\{  \text{Im } N_kW(z) \: \: \big | \: \:
k \in \{1,...,\kappa\}, \: W \in \Gamma \Bigr\}.
\end{equation}
(Cf. \cite[(3.8)]{IW}.)
This definition does not depend on the choice of $\FF$ or of the
maps $N_j$.
In fact, we have
\begin{align} \label{YGGINTRINSIC}
\ygg (z) 
= \sup  \Bigl\{ \text{Im } NW(z) \: \: \big | \: \: 
N \in \minss, \: W \in \Gamma \Bigr\},
\end{align}
where $\minss$ is the set of all $N \in \PSL(2,\R)$ such that
$\eta = N^{-1} \infty$ is a cusp and 
$\Gamma_{\eta} = N^{-1} \left [ \matr 1101 \right ] N.$
Relation \eqref{YGGINTRINSIC} follows directly from the fact that
each $N \in \minss$ has a factorization $N=\matr 1x01 N_j A$,
for some $x \in \R$, $A \in \Gamma$.

The function $\ygg(z)$ 
is well-known to be continuous and $\Gamma$-invariant,
and bounded from below by a positive constant which only depends on
the group $\Gamma$.
Notice that we have $\ygg(z) \to \infty$ as $z \in \FF$ 
approaches any of the cusps.
We also remark that, for any constant $0 \leqq c<1$:
\begin{equation} \label{YGGINTEGRAL}
\int_{\F} \ygg(z)^c \,  d\mu(z) < \infty.
\end{equation}
This is easily seen by splitting the region $\F$ into $\kappa$
cuspidal regions 
$\CC_k = N_k^{-1} \bigl( [0,1] \times [B, \infty) \bigr)$
with $B$ large$\emdash$cf.\ \eqref{CKRAV}$\emdash$
and a remaining compact region,
and then using the fact that $\ygg(z) = \tim N_k(z)$ for all 
$z \in \CC_k$. 

When proving Theorem \ref{MAINTHM} 
we will first assume $f \in C^{\infty}(T_1\HH) \cap C_c(\Tqspace)$,
i.e., that $f$ is infinitely differentiable as a function on $T_1\HH$
and has compact support as a function on $\Tqspace$.
Following \cite[p.\ 725]{Sarnak}, we start by
applying Fourier expansion in the variable $\theta$:
\begin{align} \label{THFOURIER}
& f(z,\theta) = \sum_{v \in \Z} \widehat{f}_v(z) e^{\ii v \theta},
\\ \label{THFOURIERCOEFF}
& \text{where} \quad
\widehat{f}_v(z) = \frac 1{2\pi} \int_0^{2 \pi} f(z,\theta) e^{-\ii v\theta} 
\, d\theta.
\end{align}
A simple computation then shows that
\begin{equation} \label{WEIGHTK}
\widehat{f}_v \left( \frac{az + b}{cz + d} \right) = 
\frac{(cz + d)^{2v}}{|cz+d|^{2v}} \: \widehat{f}_v(z),
\qquad
\text{for all } \matr abcd \in \Gamma, \: z \in \HH.
\end{equation}
We will call any function on $\HH$ satisfying the automorphy relation 
in \eqref{WEIGHTK} \textit{a function of weight $2v$ on $\qspace$},
and we will use $C(\qspace, 2v)$, $C_c(\qspace, 2v)$,
$\Ll_2(\qspace, 2v)$ etc.\ to denote the corresponding function 
spaces. 
Hence in our case we have 
\begin{equation*}
\widehat{f}_v \in C^{\infty}(\HH) \cap C_c(\qspace, 2v).
\end{equation*}

Let us fix $v \in \Z$ temporarily.
The function $\widehat{f}_v$ has a spectral expansion
with respect to the ``weight-$2v$-Laplacian'',
\begin{equation*}
\Delta_{2v} = y^2 \left( \frac{\partial^2}{\partial x^2} + 
\frac{\partial^2}{\partial y^2} \right) 
- 2\ii vy \frac{\partial}{\partial x}.
\end{equation*}
This operator $\Delta_{2v}$ acts on functions of weight $2v$ on
$\qspace$.
We let $\phi_0, \phi_1, ...$ be the discrete eigenfunctions of
$-\Delta_{2v}$, taken to be orthonormal
and to have increasing eigenvalues 
$\lambda_0 \leqq \lambda_1 \leqq ...$.
Cf.\ \cite[p.\ 370(items 1-3)]{TR2};
in particular, one knows that $|v| (1-|v|) \leqq \lambda_0$.
We let $\Ll_2(\lambda, 2v)$ denote the subspace of
$\Ll_2(\qspace, 2v)$ generated by the $\phi_j$'s with 
$\lambda_j = \lambda$.
The functions in $\Ll_2(\lambda, 2v)$ are usually called
\textit{Maass waveforms of weight $2v$,}
and they all belong to $C^{\infty}(\HH) \cap \Ll_2(\qspace, 2v)$.
According to \cite[pp.\ 317(Prop.\ 5.3), 414 (lines 12-16)]{TR2},
we now have the following spectral expansion,
for any given $f_v \in C^{\infty}(\HH) \cap C_c(\qspace, 2v)$:
\begin{equation} \label{SPEKTRALUTV}
f_v(z) =
\sum_{n \geqq 0} d_n \phi_n(z)
+ \sum_{k=1}^{\kappa} \int_0^{\infty} g_k(t) E_k ( z , \sfrac 12
+ \ii t, 2v) \, dt,
\end{equation}
with uniform and absolute convergence over
$z \in \HH$-compacta.
Here $E_k(z,s,2v)$ is the Eisenstein series of weight $2v$ associated
to the cusp $\eta_k$, 
cf.\ \cite[pp.\ 355(Def.\ 5.3), 368(5.19)]{TR2}.
The coefficients $d_n$ and $g_k(t)$ are given by
$d_n = \langle f_v, \phi_n \rangle$
and $g_k(t) = \frac 1{2\pi} 
\int_{\F} f_v(z) \overline{E_k(z,\frac 12+\ii t,2v)} \, d\mu (z)$
(cf.\ \cite[p.\ 243(Remark 2.4)]{TR2}).
In particular, each $g_k(t)$ is a continuous function on $[0,\infty)$.

The proof in \cite{TR2}
of the beautiful convergence in \eqref{SPEKTRALUTV}  
starts by considering the spectral expansion (in $L_2$-sense) 
of the function $\Delta_{2v} f_v + 
a(1-a)f_v \in L_2(\qspace,2v)$,
for some fixed number $a > |v|+1$; 
this is then integrated against the Green's function $G_a(z,w,2v)$.
It is seen in this proof that
\begin{align} \label{SPEKTRALBOUND} 
\sum_{n \geqq 0}
& | d_n |^2 \left( a(1-a) - \lambda_n \right) ^2
+ 2\pi \sum_{k=1}^{\kappa} \int_0^{\infty} 
\left | g_k(t) \right |^2 \left( a(1-a) - \sfrac 14 - t^2 \right) ^2 
\, dt
\\ \notag
& = \int_{\F} \left | \Delta_{2v} f_v(z) + a(1-a) 
f_v(z) \right | ^2 <
\infty. 
\end{align}
Cf. \cite[pp.\ 91(9.36), 244--245]{TR2}.

We will need a bound on the rate of convergence in  
\eqref{SPEKTRALUTV} which is  \textit{uniform} over all $z \in \HH$.
This is obtained in the next two lemmas,
the first of which is a generalization of \cite[Prop.\ 7.2]{IW}.

\begin{lemma} \label{BESSELINEQ}
Given $v \in \Z$ and
$\phi_0, \phi_1,...$ as above,
we have, for all $z \in \HH$ and $T \geqq 1$: 
\begin{align} \label{BESSELINEQRES}
\sum_{\lambda_n \leqq \frac 14 + T^2} | \phi_n(z)|^2
& + \sum_{k=1}^{\kappa} \int_0^T
|E_k(z, \sfrac 12 + \ii t, 2v)|^2 \, dt
\\ \notag
& = O \Bigl( (T+|v|) \ygg(z) + (T+|v|)^2 \Bigr).
\end{align}
The implied constant depends only on $\Gamma$,
and \textit{not} on $v,T,z$.
\end{lemma}

The uniformity in $v$ in the above bound will not be essential in 
the proofs of the main results in the present paper.

\begin{Proof}
For $v=0$, this is Proposition 7.2 in \cite{IW}.
We will now assume $v \in \Z$, $v \neq 0$,
and will show how to carry over the proof in \cite{IW} to this case.

We let $\chi_{\delta}$ be the characteristic function of the interval 
$[0,\delta]$, where $\delta$ is a (small) positive constant
to be specified later.
We define, for $z,w \in \HH$,
\begin{align*}
& u(z,w) = \frac{|z-w|^2}{4 \tim z \tim w}
\qquad \text{(cf.\ \cite[p.\ 8(1.4)]{IW})};
\\
& k(z,w) = 
(-1)^v \frac{(w-\overline{z})^{2v}}{|w-\overline{z}|^{2v}}
\: \: \chi_{\delta} \bigl( u(z,w) \bigr) ;
\\
& K(z,w) = \sum_{T \in \Gamma} k(z,Tw) \:
\frac{(cw+d)^{2v}}{|cw+d|^{2v}}
\quad (\text{wherein } \: T= \matr **cd).
\end{align*}
This agrees with \cite[pp.\ 359--360, (2.7)]{TR1}
(for trivial character and ``$\Phi = \chi_{4 \delta}$'').
As in \cite[pp.\ 107--109]{IW}, 
we fix $w \in \HH$ and consider $K(z,w)$ as a function of $z$;
this function belongs to $L_2(\qspace,2v)$
and has compact support in $\qspace$.
We may now apply Bessel's inequality to obtain
\begin{align} \label{BESSELINEQ2}
& \sum_{n \geqq 0} |h(t_n) \phi_n(w)|^2
+ \frac 1{2\pi}  \sum_{k=1}^{\kappa} \int_0^{\infty} 
\bigl | h(t) E_k(w,\frh+\ii t,2v)  \bigr | ^2 \, dt
\\ \notag
& \leqq
\int_{\F} |K(z,w)|^2 \, d\mu(z).
\end{align}
Here
\begin{equation} \label{BESSELINEQ3}
h(t) = \int_{\HH} k(\ii,z) y^{\frac 12 + \ii t} \, d\mu(z),
\end{equation}
and the $t_n$'s are defined through
$\lambda_n = s_n(1-s_n)$, $s_n = \frac 12 + \ii t_n$,
with $s_n \in \bigl( \frac 12 + \ii[0,\infty) \bigr)
\cup (\frac 12, |v|]$
(recall here that $\lambda_n \geqq |v|(1-|v|)$).
The proof of \eqref{BESSELINEQ2} and \eqref{BESSELINEQ3}
uses \cite[pp.\ 291(3.23), 373(item 12)]{TR2}
and unfolding of the integral
$\langle K(\cdot , w), \phi_n \rangle = 
\overline{\int_{\qspace} K(w,z) \phi_n(z) \, d\mu(z)}$
and the analogous Eisenstein integral,
together with an application of 
\cite[p. 364(Prop.\ 2.14)]{TR1}.
Notice here that the proof in
\cite{TR1} of Prop.\ 2.14 remains valid if the assumption in
\cite[p.\ 359(Def.\ 2.10)]{TR1} is replaced by the weaker assumption that
$\Phi$ is piecewise continuous and of compact support.

By an argument exactly as in \cite[pp.\ 109--110]{IW} we have
the following upper bound:
\begin{equation} \label{BESSELINEQ1}
\int_{\F} |K(z,w)|^2 \, d\mu(z) = O(\delta^{\frac 32} \ygg(w) +
\delta).
\end{equation}

Next, we want to bound $h(t)$ from below. 
Let us assume
\begin{equation} \label{DELTAASS}
\delta \leqq 
(100|v|)^{-2} 
\qquad
\text{(hence in particular, $\delta \leqq 10^{-4}$).}
\end{equation}
Let $D$ be the hyperbolic disc defined by $u(\ii,z) \leqq \delta$.
One easily checks that $|z-\ii| < 3 \sqrt{\delta}$ holds for all
$z \in D$, and hence
\begin{align*}
\bigl | \Arg(z+\ii) - \sfrac{\pi}2 \bigr |
\leqq \arcsin( 3\sqrt{\delta} /2) \leqq 3 \sqrt{\delta}.
\end{align*}
Hence, for all $z \in D$, we have
$\bigl | \Arg \: \, k(\ii,z) \bigr | \leqq 6|v| \sqrt{\delta} <
\sfrac{\pi}3$
(cf.\ \eqref{DELTAASS}),
and thus $\tre k(\ii,z) \geqq \frac 12$,
since $|k(\ii,z)|=1$.
Using the fact that the hyperbolic area of $D$ is
$4\pi \delta$, we obtain
\begin{equation*}
\left | h( \sfrac i2 ) \right | 
\geqq \int_{\HH} \tre k(\ii,z) \, d\mu(z) \geqq 
2 \pi \delta.
\end{equation*}
Now take any $t \in \C$ such that
$s = \frac 12 + \ii t \in
\bigl( \frac 12 + \ii[0,\infty) \bigr) \cup (\frac 12, |v|]$.
Notice that $z = x+{\ii y} \in D$ implies
$|y-1| < 3 \sqrt{\delta} < |2v|^{-1}$,
and in this $y$-interval we have
\begin{align*}
\left | \frac d{dy} (y^{s}) \right |
& = |s| y^{\tre s-1} 
\leqq |s| \max \left( (1-|2v|^{-1})^{-\frac 12}, (1+|2v|^{-1})^{|v|-1}
\right) 
\\
& < |s| \max \left( 2, \exp \bigr\{ (|v|-1) |2v|^{-1} 
\bigr\} \right) < 3|s|.
\end{align*}
(since $\log(1+|2v|^{-1}) < |2v|^{-1}$).
It follows that 
$\bigl | y^s - 1 \bigr | \leqq 3 |s| |y-1| 
\leqq 9 \sqrt{\delta} |s|$ for all $z \in D$.
We now keep $T \geqq |v|$, and let
$\delta = (100T)^{-2}$
(notice that \eqref{DELTAASS} is then fulfilled).
We then have $9 \sqrt{\delta} |s| < \frac 14$ 
for all $s = \frac 12 + \ii t \in \frac 12 + \ii[0,T]$,
and also for all $s \in [\frac 12,|v|]$.
Hence, for all these $s$, we have by \eqref{BESSELINEQ3}:
$| h(t) - h(\sfrac i2) | \leqq \sfrac 14 \mu(D)
\leqq \pi \delta,$ and thus
\begin{equation*} 
|h(t)| \geqq \pi \delta.
\end{equation*}
Combined with \eqref{BESSELINEQ2} and \eqref{BESSELINEQ1}, this gives:
\begin{align*}
\sum_{\lambda_n \leqq \frac 14 + T^2} | \phi_n(w)|^2
& + \sum_{k=1}^{\kappa} \int_0^T
|E_k(w, \sfrac 12 + \ii t, 2v)|^2 \, dt
\\
& = 
O(\delta^{-\frac 12} \ygg(w) + \delta^{-1})=
O \Bigl( T \ygg(w) + T^2 \Bigr).
\end{align*}
This holds for all $T \geqq |v|$.
The desired inequality \eqref{BESSELINEQRES} now follows,
using the fact that the left hand side in
\eqref{BESSELINEQRES} is an increasing function of $T$.
\end{Proof}

\begin{lemma} \label{HUGGSPEKTRAL}
Let $v \in \Z$ and $f_v \in C^{\infty}(\HH) \cap C_c(\qspace, 2v)$ be
given. We then have in \eqref{SPEKTRALUTV},
for all $T \geqq1$, $z \in \HH$:
\begin{align} \label{HUGGSPEKTRAL1}
& \left | f_v(z) 
- \sum_{\lambda_n \leqq 1/4 + T^2} d_n \phi_n(z)
-  \sum_{k=1}^{\kappa} \int_0^T g_k(t) E_k(z,\sfrac 12 + \ii t,2v) \, dt
\right |  
\\ \notag
& \leqq
O \left( T^{-1} +  T^{-3/2} \sqrt{\ygg(z)} \right).
\end{align}
The implied constant depends on $\Gamma$, $v$ and $f_v,$ 
but {\rm not} on $T,$ $z$.
\end{lemma}

\begin{Proof}
By \eqref{SPEKTRALUTV}, the difference in the left hand side 
of \eqref{HUGGSPEKTRAL1} equals
\begin{equation} \label{HUGGPROOF1}
\sum_{\lambda_n > 1/4 + T^2} d_n \phi_n(z)
+ \sum_{k=1}^{\kappa} \int_T^{\infty} g_k(t) E_k(z,\sfrac 12 + \ii t,2v) \,
dt.
\end{equation}
By \eqref{SPEKTRALBOUND}, the sum and the $\kappa$ integrals
\begin{equation*}
\sum_{\lambda_n > 1/4 + 1} | d_n |^2 \lambda_n^2;
\qquad
\int_1^{\infty} |g_k(t)|^2 t^4 \, dt 
\quad (k=1,2,...,\kappa)
\end{equation*}
are all bounded from above by some finite constant 
(which depends on $v$, $f_v,$ $\Gamma$).
Hence by writing $d_n \phi_n(z) = \left( d_n \lambda_n \right) \cdot
\left( \phi_n(z) / \lambda_n \right)$ 
and $g_k(t) E_k(...) = \left( g_k(t) t^2 \right) \cdot 
\left( E_k(...) / t^2 \right)$ in \eqref{HUGGPROOF1},
and then applying Cauchy's inequality, we find that the modulus of
\eqref{HUGGPROOF1} is bounded from above by
\begin{align} \notag
& O(1) \sqrt{\sum_{\lambda_n > 1/4 + T^2}
\lambda_n^{-2} |\phi_n(z)|^2}
+ O(1) \sum_{k=1}^{\kappa} \sqrt{ 
\int_T^{\infty} t^{-4} |E_k(z,\sfrac 12 + \ii t,2v)|^2 \, dt}
\\ \notag
& = O(1) \sqrt{
\sum_{t_n > T}
t_n^{-4} |\phi_n(z)|^2
+ \sum_{k=1}^{\kappa}
\int_T^{\infty} t^{-4} |E_k(z,\sfrac 12 + \ii t,2v)|^2 \, dt},
\end{align}
where the $t_n$'s are as in the proof of Lemma \ref{BESSELINEQ}. 
Now \eqref{HUGGSPEKTRAL1} follows from Lemma \ref{BESSELINEQ}
(applied for our \textit{fixed} $v$, so that the right hand side in 
\eqref{BESSELINEQRES} is $O(T \ygg(z) + T^2)$),
and partial summation.
\end{Proof}

Next, we will review some facts about the Fourier expansion of Maass
waveforms of even integer weight.
We will use the standard notation $W_{k,m}(z)$ 
for the Whittaker function as in (e.g.) 
\cite[Ch.\ 7 (10.04),(11.03)]{Olver}.
$W_{k,m}(z)$ is a holomorphic function for 
$k,m \in \C$, 
$|\text{Arg}(z)|< \pi$.

\begin{lemma} \label{EGENVLEMMA}
Let $\phi \in \Ll_2(\lambda, 2v)$, $\phi \not\equiv 0$,
$v \in \Z$;
we can then write $\lambda = s(1-s)$ for a unique
$s \in \frac 12 + \ii[0, \infty)$ or $s \in (\frac 12, \max(1,|v|)]$.
We have a Fourier expansion
\begin{equation} \label{VFEXP}
\phi(z) = 
c_0 y^{1-s} +
\sum_{n \neq 0} \frac {c_n}{\sqrt{|n|}} W_{v \cdot \sgn(n) ,s-\frac 12} 
(4 \pi |n| y) e(nx).
\end{equation}
From this it follows that, as $y \to \infty$,
\begin{equation} \label{EGENVLEMMA2}
\phi(z) = 
c_0 y^{1-s} + O(y^{|v|} e^{-2\pi y} ).
\end{equation}
If $\tre s = \frac 12$, then 
we necessarily have $c_0=0$.
\end{lemma}

\begin{Proof}
This is a restatement of 
\cite[pp.\ 370-371(items 1,3,5)]{TR2}
and \cite[pp.\ 26(Prop 4.12), 349]{TR2}.
Concerning the translation of \cite[p.\ 370(5.28)]{TR2}
into \eqref{VFEXP},
cf.\ \cite[pp.\ 347(Lemma 4.5), 348(lines 3-5), 355(Def.\ 5.4)]{TR2}
and \cite[p.\ 256 (footnote $\ddagger$)]{Olver}.
\end{Proof}

We remark that the Whittaker function is always
real valued in \eqref{VFEXP}; 
in fact, we have $W_{k,m}(z) \in \R$ for all
$k \in \R$, $m \in \R \cup i \R$, $z>0$,
as follows from \cite[Ch.\ 7 (9.03),(11.02),Ex.\ 11.1]{Olver}. 

Each Maass waveform of even integer weight
allows an explicit description in terms  
of either a Maass waveform of weight zero or a holomorphic
cusp form, and the coefficients in the Fourier expansions of the two
forms are proportional.
We will state this fact in a precise form in the next two lemmas,
thus slightly generalizing formulas given earlier in \cite{DJ94,DJ97}.
For the proofs of the main results in the present paper
it will not be essential to know the exact formulas for the
proportionality factors involved.
However, these formulas, as well as the uniformity in $v$ obtained
in Lemma \ref{BESSELINEQ}, might be of importance in future
applications, e.g.\ to obtain results on \textit{the rate of convergence} 
in various asymptotic geometric problems.
In this vein, note that the (uniform) asymptotic properties
of $W_{k,m}(z)$ are known for all relevant ranges of $k,m,z$;
cf.\ \cite{OlverW,Dunster1,Dunster2}.

\enlargethispage{25pt}

We use the following standard notation: 
\begin{align*}
& \KKK_{2v} = iy \frac{\partial}{\partial x} 
+ y \frac{\partial}{\partial y} + v,
\qquad
\LLL_{2v} = iy \frac{\partial}{\partial x} 
- y \frac{\partial}{\partial y} + v,
\\
& (\alpha)_n = \Gamma(\alpha+n)/\Gamma(\alpha) = 
\alpha (\alpha+1) ... (\alpha+n-1)
\qquad
\text{for } \: \alpha \in \C, \: n \geqq 0.
\end{align*}

\begin{lemma} \label{MAASSWV}
Let $\phi$ be as in Lemma \ref{EGENVLEMMA}, 
and assume 
$\tre s < 1$. 
Then there exists a Maass waveform $\phi_0$ of weight $0$ on $\qspace$ 
satisfying $\Delta_0 \phi_0 + \lambda \phi_0 = 0$ and
$|| \phi_0 ||_{L_2} = || \phi ||_{L_2}$,
such that
\begin{align*}
\phi = 
\bigl\{ (s)_{|v|} (1-s)_{|v|} \bigr\}^{-\frac 12}
\begin{cases}
\KKK_{2v-2} \KKK_{2v-4} ... \KKK_0 (\phi_0) 
& \text{if } v \geqq 1
\\
\phi_0 & \text{if } v = 0
\\
\LLL_{2v+2} \LLL_{2v+4} ... \LLL_0 (\phi_0)
& \text{if } v \leqq -1.
\end{cases}
\end{align*}
Let us, furthermore, 
assume $\phi_0$ to have the following
Fourier expansion
in the standard format involving the Macdonald-Bessel function
$K_{\mu}(z)$:
\begin{equation*}
\phi_0(z) = d_0 y^{1-s} +
\sum_{n \neq 0} d_n \sqrt y K_{s-\frac 12}(2 \pi |n| y) e(nx)
\end{equation*}
(cf.\ \cite[Ch.\ 6, \S 4]{TR2}).
We then have in \eqref{VFEXP}:
\begin{align} \label{COEFF0REL}
& c_0 = \begin{cases}
\sgn(v)^v \bigl\{ (1-s)_{|v|} / (s)_{|v|} \bigr\}^{\frac 12} \, d_0
& \text{if } v \neq 0 \text{ and } \frac 12 < s < 1
\\ 
d_0 & \text{otherwise} \end{cases}
\\ \notag
& c_n = \frac 12
\bigl\{ -\sgn(n) \bigr\}^v 
\bigl\{ (s)_{|v|} (1-s)_{|v|} \bigr\}^{- \frac 12 \sgn(nv)}  
\, d_n \quad \text{for } \:  n \neq 0.
\end{align}
\end{lemma}

Notice that all square roots appearing in the lemma are
well defined, since clearly $(s)_{|v|} (1-s)_{|v|} > 0$
for all 
$s \in \frac 12 + \ii[0, \infty)$ and for all $s \in (\frac 12, 1)$,
and 
$(1-s)_{|v|} / (s)_{|v|} > 0$ for all $s \in (\frac 12,1)$.
We also remark that the above formulas are consistent with 
\cite[(1.8)]{DJ94}, \cite[(6),(7)]{DJ97}
in the case $\tre s = \frac 12$, $\Gamma = \text{PSL}(2,\Z)$ 
(since the $\varphi_{j,k}$-expansions in \cite{DJ94,DJ97} agree with
our formulas above after multiplication by a certain constant
$c = c(s,v) \in \C$ with $|c|=1$).

\begin{Proof}
The first assertion follows from
\cite[p.\ 382(d),(e),(f),(h)]{TR2},
used repeatedly.

It remains to prove the formulas for the $c_n$'s.
The case $v=0$ is trivial since
$K_{s-\frac 12}(y) = (\pi / 2y)^{\frac 12} W_{0,s-\frac 12}(2y)$
(cf.,\ e.g.,\ \cite[Ch.\ 7 Ex.\ 10.1,(11.03)]{Olver}).
The cases $v>0$ and $v<0$ are now treated by first
noticing that we may differentiate term-by-term any
number of times in \eqref{VFEXP}
(cf.\ \cite[pp.\ 25(Remark 4.11), 349(line 1)]{TR2}),
and then making repeated use of 
\cite[Lemma 1.1]{DJ94}. (Notice here that \cite[Lemma 1.1]{DJ94}
can be extended to \textit{complex} $t$.
Notice also that the right hand side in
the second formula in \cite[Lemma 1.1]{DJ94} 
should be corrected 
by interchanging ``$k-1$'' and ``$-(k-1)$'' in the arguments of the
Whittaker functions.)
%
\end{Proof}

Before stating the next lemma,
we recall that if $m \in \{2,4,6,...\}$, and
if $\phi_0$ is a holomorphic cusp form 
of weight $m$ on
$\qspace$, then $y^{\frac m2} \phi_0$ and $y^{\frac m2} 
\overline{\phi_0}$
are Maass waveforms of weight $m$ and $-m$:
\begin{equation*}
y^{\frac m2} \phi_0 \in L_2(\lambda,m),
\quad y^{\frac m2} \overline{\phi_0} \in L_2(\lambda,-m),
\quad \text{where } \: \lambda = \sfrac m2 (1- \sfrac m2).
\end{equation*}
In particular, 
$|| y^{\frac m2}  \phi_0 ||_{L_2} 
= || y^{\frac m2} \overline{\phi_0} ||_{L_2} = 
\sqrt{\int_{\F} |\phi_0(z)|^2 y^{m} \, d\mu(z)}$.
Cf.\ \cite[p.\ 382(item 22), Prop 5.14]{TR2}.

\begin{lemma} \label{HOLO}
Let $\phi$ be as in Lemma \ref{EGENVLEMMA}, 
and assume 
$s \geqq 1$. 
Then, if $v=0$, we must have $s=1$, and $\phi$ is a constant function.
Now assume $v \neq 0$.
We then 
have $s \in \{1,2,...,|v| \}$, and
there exists a holomorphic cusp form $\phi_0$ on $\qspace$ of weight
$2s$ with $|| y^s \phi_0 ||_{L_2} = || \phi ||_{L_2}$,
such that
\begin{align*}
\phi = 
\bigl\{ (|v|-s)! (2s)_{|v|-s} \bigr\}^{-\frac 12}
\begin{cases}
\KKK_{2v-2} \KKK_{2v-4} ... \KKK_{2s} \bigl( y^{s} \phi_0 \bigr)
& \text{if } v \geqq 1
\\
\LLL_{2v+2} \LLL_{2v+4} ... \LLL_{-2s} 
\Bigl( y^{s} \overline{\phi_0} \Bigr) 
& \text{if } v \leqq -1.
\end{cases}
\end{align*}
(If $s = |v|$, this should be interpreted as
$\phi = y^s \phi_0$ if $v \geqq 1$, $\phi = y^s \overline{\phi_0}$ if
$v \leqq -1$.)
Furthermore, if $\phi_0$ has Fourier expansion
\begin{equation*}
\phi_0(z) = \sum_{n=1}^{\infty} d_n e(nz),
\end{equation*}
then we have in \eqref{VFEXP}:
\begin{align*}
& c_n = \bigl\{ (|v|-s)!(2s)_{|v|-s} \bigr\} ^{-\frac 12}
(4 \pi)^{-s} |n|^{\frac 12 - s} \,
\begin{cases}
(-1)^{v-s} d_n & \text{if } v>0
\\
\overline{d_{-n}} & \text{if } v<0
\end{cases}
\end{align*}
for all $n$ with $nv > 0$;
also,
$c_0 = 0$, and $c_n = 0$ for all $n$ with $nv < 0$.
\end{lemma}

\begin{Proof}
The case $v=0$ is trivial, cf.\ \cite[p.\ 71(Claim 9.2)]{TR2}.
If $v \neq 0$, then we have $s \in \{1,2,...,|v|\}$
by \cite[pp.\ 350(5.4), 427(line 1)]{TR2}.
The rest of the proof is very similar to the proof of Lemma
\ref{MAASSWV},
except that we also use \cite[p.\ 383(a),(d),(g)]{TR2},
and the fact that
\begin{equation} \label{WHITTELEM}
W_{s,s-\frac 12}(y) = e^{-y/2} y^s,
\end{equation}
which follows easily from (e.g.) 
\cite[Ch.\ 7 (9.03),(9.04),(10.09),(11.03)]{Olver}.
\end{Proof}

The above Lemmas \ref{EGENVLEMMA},
\ref{MAASSWV}, \ref{HOLO} 
deal with the Fourier expansion of $\phi(z)$ at the cusp
$\eta_1=\infty$.
Of course, there are analogous results for
the Fourier expansions corresponding to
the other cusps $\eta_2,...,\eta_{\kappa}$ of $\Gamma$
(these results may be proved e.g.\ by applying the above lemmas to the
Fuchsian groups $N_k \Gamma N_k^{-1}$, $k=2,...,\kappa$).
One consequence of this is the following:
For any $\phi \in L_2(\lambda,2v)$ 
with $\lambda = s(1-s)$ as in Lemma \ref{EGENVLEMMA}, 
there exists a constant $C>0$ such that
\begin{equation} \label{MWGROWTH}
|\phi(z)| \leqq C \ygg(z)^{1 - \tre s},
\qquad \forall z \in \HH.
\end{equation}
(For notice that for each $k \in \{1,...,\kappa\}$,
the $\eta_k$-analog of \eqref{EGENVLEMMA2} 
implies that $\phi(z) = O((\text{Im } N_kz)^{1-\tre s})$ 
as $z \to \eta_k$ inside $\FF$.
Hence there is a $C>0$ such that \eqref{MWGROWTH} holds for all 
$z \in \FF$. By $\Gamma$-invariance, \eqref{MWGROWTH} now holds
for all $z \in \HH$.)

\vspace{5pt}

One main ingredient in our proof of Theorem \ref{MAINTHM} will be the
well-known Rankin-Selberg type bound on the sums 
$\sum_{|n| \leqq N} |c_n|^2$ of the Fourier coefficients of the
Maass waveforms (cf.,\ e.g.,\ 
\cite[Thm.\ 3.2]{IW}, \cite[Thm.\ 5.1]{IW2}).
We will also need a similar bound on the Fourier coefficients of the
Eisenstein series, which we will prove in 
Lemma \ref{EISENSTEINRS} below.

\enlargethispage{20pt}

We first recall the following 
explicit formula for the Eisenstein series of even
integer weight in terms of the Eisenstein series of weight zero:

\begin{lemma} \label{EISENST}
Consider the Fourier expansion of the Eisenstein series
of weight zero:
\begin{align*}
E_k(z,s,0) = \delta_{k1} y^s + \varphi_{k1}(s) y^{1-s}
+ \sum_{n \neq 0} \psi_{n,k}(s) \sqrt y K_{s-\frac 12} (2\pi |n| y) 
e(nx).
\end{align*}
We then have the following expansion of the Eisenstein series of
weight $2v$:
\begin{align*}
E_k(z,s,2v) & = \delta_{k1} y^s + 
(-1)^v \frac{\Gamma(s)^2}{\Gamma(s+v) \Gamma(s-v)} \varphi_{k1}(s)
y^{1-s}
\\
& + \sum_{n \neq 0} 
\frac {(-1)^v \Gamma(s)}{2 \Gamma(s+v \cdot\sgn(n))} 
 \frac {\psi_{n,k}(s)}{\sqrt{|n|}} W_{v \cdot \sgn(n), s-\frac 12} 
(4\pi |n| y) e(nx).
\end{align*}
\end{lemma}

\begin{proof}
For $\tre s > 1$ this follows from the explicit formulas in 
\cite[pp.\ 368(5.22), 369(5.23)]{TR2}.
(Cf.\ also our remarks in the proof of Lemma \ref{EGENVLEMMA},
as well as \cite[p.\ 280(Prop 3.7)]{TR2}.)
For general $s$ the result now follows by meromorphic continuation,
since it is known that, for any fixed $k \in
\{1,...,\kappa\}$ and $v \in \Z$,
$E_k(z,s,2v)$ has a Fourier expansion of the above type, with each
coefficient being a meromorphic function in $s$.
Cf.\ \cite[p.\ 374(items 15-17)]{TR2}.
\end{proof}

\begin{lemma} \label{EISENSTEINRS}
Let $v \in \Z$, $k \in \{1,...,\kappa\}$
and $s = \sigma + it \in \C$, $\sigma \geqq \frac 12$.
We assume that the Eisenstein series $E_k(z,s,2v)$ does not
have a pole at $s$, and we take the Fourier expansion to be
\begin{align} \label{EFOURIEREXP} 
E_k(z,s,2v) & = \delta_{k1} y^s + c_0 y^{1-s}
+ \sum_{n \neq 0} 
\frac {c_n}{\sqrt{|n|}} W_{v \cdot \sgn(n), s-\frac 12} 
(4\pi |n| y) e(nx).
\end{align}
We then have, for $N \geqq 1$:
\begin{equation} \label{EISENSTEINRSRES}
\sum_{1 \leqq |n| \leqq N} |c_n|^2 = 
O(1) \begin{cases}
N \log 2N & \text{if } \sigma = \frac 12
\\
N^{2\sigma} & \text{if } \sigma > \frac 12.
\end{cases}
\end{equation}
(The implied constant depends only on $\Gamma$, $v$, $s$.)
\end{lemma}

\begin{Proof}
In the case $\sigma = \frac 12$, $v=0$,
the desired result was proved in \cite[Lemma 2.3.8]{Strhoro}.
We will review the proof from \cite{Strhoro}, and show that it
generalizes to $\sigma \geqq \frac 12$, 
$v \in \Z$.

We recall that there are only a finite number of
poles of $E_k(z,s,0)$ in the half plane $\sigma \geqq \frac 12$, 
and all these poles belong to the interval $(\frac 12, 1]$.
It now follows from the formulas in Lemma \ref{EISENST} that 
if $s \neq 1$, $\tre s \geqq \frac 12$, and $s$ is not a pole of
$E_k(z,s,2v)$, then $s$ cannot be a pole of $E_k(z,s,0)$;
these formulas also imply (for $s \neq 1$) that it suffices to prove 
\eqref{EISENSTEINRSRES} for $v=0$.
However, as we will see,
it is just as easy to prove Lemma \ref{EISENSTEINRS} directly
for general $v \in \Z$, and this also covers the case $s=1$
($v \neq 0$).

We keep $0<Y<H$ and study
the following integral:
\begin{equation} \label{EISENSTEINRSPROOF8}
J = \int_{\DD} \left | E_k ( z,s,2v ) \right |^2 \, d\mu (z),
\qquad
\text{where} \quad \DD = (0,1) \times (Y,H).
\end{equation}
We tessellate $\DD$ by translates of the fundamental region,
i.e.\ we write $\DD=$\linebreak $\cup_{T \in \Gamma} (\DD \cap T(\F)),$
an essentially disjoint union.
It then turns out that $\DD$ is fully covered by the
translates of a \textit{truncated} region
\begin{equation*}
\F_B = \F - \bigcup_{j=1}^{\kappa} 
N_j^{-1} \bigl( [0,1] \times [B, \infty ) \bigr),
\quad \text{where } \: B = \max \left( B_0,H,\frac 1Y \right).
\end{equation*}
(Cf.\ \eqref{CKRAV}. For details, see the proof in
\cite{Strhoro}.)
The upshot of this is that
\begin{equation*}
J = \int_{\F_B} \# \Bigl\{ T \in \Gamma \mid T(z) \in \DD
\Bigr\} \cdot \left| E_k (z,s,2v) \right | ^2
\, d\mu (z).
\end{equation*}
But we have by \cite[Lemma 2.10]{IW}:
\begin{align*}
\# \Bigl\{ T \in \Gamma \quad \Big | \quad T(z) \in \DD \Bigr\}
& \leqq \# \Bigl\{ W_0 \in [S] \setminus \Gamma \quad \Big | \quad 
\tim W_0(z) > Y \Bigr\} 
\\
& = 1+ O \left( Y^{-1} \right) ,
\end{align*}
where the implied constant depends only on $\Gamma$,
i.e.\ the bound is uniform over all $z \in \HH$ and $Y>0$.
Hence,
\begin{equation} \label{EISENSTEINRSPROOF5}
J = O\bigl( 1+ Y^{-1} \bigr) \int_{\F_B} 
\left| E_k (z,s,2v) \right | ^2
\, d\mu (z).
\end{equation}
We now decompose $\F_B$ as a union of the compact region
$\F_{B_0}$ and the cuspidal regions
$N_j^{-1} \bigl( [0,1] \times [B_0, B ) \bigr)$,
for $j=1,...,\kappa$.
From the Fourier expansion of $E_k(z,s,2v)$ in the cusp $\eta_j$ it
follows that for all points $z = N_j^{-1}(x'+y'i)$, $y' \geqq B_0$ 
we have $E_k(z,s,2v) = O((y')^{\sigma})$
(cf.\ 
Lemma \ref{EISENST} and the proof of \eqref{MWGROWTH}); 
of course, the implied constant depends on $s$ and $v$.
Using this we easily obtain
\begin{equation} \label{EISENSTEINRSPROOF3}
J = O\bigl( 1+ Y^{-1} \bigr) 
\begin{cases}
\log (2B) & \text{if } \: \sigma = \frac 12
\\
B^{2\sigma-1} & \text{if } \: \sigma > \frac 12.
\end{cases}
\end{equation}

On the other hand,
substituting \eqref{EFOURIEREXP} 
directly in the definition of $J$, \eqref{EISENSTEINRSPROOF8},
and then applying Parseval's formula, we get:
\begin{align*} 
J \geqq 4 \pi \sum_{n \neq 0} |c_n|^2 
\int_{4 \pi |n| Y}^{4 \pi |n| H}
W_{v \cdot \sgn(n), s-\frac 12} (u)^2 \, \frac{du}{u^2}.
\end{align*}

We now take $H=1$ and $Y=(4 \pi N)^{-1}$.
With this choice we have
$[1,4\pi] \subset [4 \pi |n| Y,4 \pi |n| H]$
whenever $1 \leqq |n| \leqq N$,
and hence the last inequality implies
\begin{align*}
\sum_{1 \leqq |n| \leqq N} |c_n|^2 
\leqq O(J). 
\end{align*}
Using this fact and \eqref{EISENSTEINRSPROOF3},
we obtain \eqref{EISENSTEINRSRES}.
\end{Proof}

\section{Averages of eigenfunctions for generic $\alpha$}

In this section we will study the sum $\frac 1M \sum_{m=1}^M f(m \alpha + iy)$ in  
the case when $f$ is an eigenfunction on $\qspace$, i.e., 
either a Maass waveform or an Eisenstein series.

\begin{proposition} \label{SMAASSWV}
Let $\varepsilon>0$, and
let $\phi$ be a non-constant Maass waveform of even integer weight,
say $\phi \in L_2(\lambda,2v)$ and $\lambda = s(1-s)$
as in \mbox{Lemma \ref{EGENVLEMMA}}.
We write
\begin{equation*}
S = S(M,y, \alpha) = \frac 1M \sum_{m=1}^M \phi(m \alpha + {\ii y}).
\end{equation*}
We then have, for all $M \geqq 1$ and $0 < y \leqq 1$,
\begin{equation} \label{SMAASSWVRES}
\int_0^1 |S|^2 \, d\alpha = O \left( y^{-\varepsilon} M^{\varepsilon -
1} \right) + \bigl[  {\textstyle \text{if $s \in (\frh , 1)$: }} \: 
 O ( y^{2-2s} ) \bigr].
\end{equation}
(The implied constant depends on $\Gamma, \phi, \varepsilon$,
but not on $M, y$.)
\end{proposition}

\begin{Proof}
As in Lemma \ref{EGENVLEMMA} we have a Fourier expansion
\begin{equation} \label{FEXPNY}
\phi(z) = 
c_0 y^{1-s} +
\sum_{n \neq 0} \frac {c_n}{\sqrt{|n|}} 
W_{v \cdot \sgn(n) ,s - \frac 12} 
(4 \pi |n| y) e(nx).
\end{equation}
It follows that
\begin{align*}
S & = 
c_0 y^{1-s} +
\frac 1M \sum_{n \neq 0} \sum_{m=1}^M 
\frac {c_n}{\sqrt{|n|}} W_{v \cdot \sgn(n) ,s - \frac 12} 
(4 \pi |n| y) e(nm\alpha)
\\
& = c_0 y^{1-s} +
\frac 1M \sum_{l \neq 0} 
\Bigl( \sum_{\substack{m \mid l \\ 1 \leqq m \leqq M}} 
\frac {c_{l/m}}{\sqrt{|l/m|}} W_{v \cdot \sgn(l) ,s - \frac 12} 
(4 \pi |l/m| y) \Bigr) e(l\alpha).
\end{align*}
We will write $s = \sigma + \ii t$,
and use $W(Y)$ as a shorthand for 
$W_{v \cdot \sgn(l) ,s - \frac 12} (Y)$.
Applying Parseval's formula and
then Cauchy's inequality, we get
\begin{align*}
\int_0^1 |S|^2 \, d\alpha & =
|c_0|^2 y^{2-2\sigma} +
\frac 1{M^2} \sum_{l \neq 0} 
\Big | \sum_{\substack{m \mid l \\ 1 \leqq m \leqq M}} 
\frac {c_{l/m}}{\sqrt{|l/m|}} W
(4 \pi |l/m| y) \Bigr |^2
\\
& \leqq
|c_0|^2 y^{2-2\sigma} +
\frac 1{M^2} \sum_{l \neq 0} \tau(|l|)
\sum_{\substack{m \mid l \\ 1 \leqq m \leqq M}} 
\frac {|c_{l/m}|^2}{|l/m|} W
(4 \pi |l/m| y)^2
\\
& \text{(where $\tau(|l|)$ is the number of divisors of $|l|$)}
\\
& = 
|c_0|^2 y^{2-2\sigma} +
\frac 1{M^2} \sum_{n \neq 0} 
\frac {|c_{n}|^2}{|n|} W_{v \cdot \sgn(n), s - \frac 12} 
(4 \pi |n| y)^2
\cdot \sum_{m=1}^M \tau(|mn|).
\end{align*}
In this sum we will use the following crude bound,
which applies whenever $n \neq 0$ and $c_n \neq 0$:
\begin{equation} \label{SMAASSWV10}
W_{v \cdot \sgn(n),s-\frac 12}(Y) = O \bigl( e^{-Y/4} \bigr),
\qquad
\forall Y > 0.
\end{equation}
(The implied constant depends on $v,s$.) 
To prove \eqref{SMAASSWV10}, first recall that
$W_{k,m}(Y) \sim e^{-Y/2}Y^{k}$ as $Y \to \infty$
(cf.\ \cite[Ch.\ 7 (11.05)]{Olver}).
It now remains to treat the case $0 < Y \leqq 1$.
If $s \in \frac 12 + \ii[0, \infty)$ or $s \in (\frac 12,1)$,
then we have, as $Y \to 0^+$:
\begin{align} \label{WFAKTA}
& W_{v \cdot \sgn(n), s - \frac 12} (Y) 
= O(1) 
\begin{cases}
Y^{\frac 12} \log (1/Y) & \text{if } s = \frac 12
\\
Y^{1-\sigma} & \text{otherwise}.
\end{cases} 
\end{align}
Cf.\ \cite[Ch.\ 7 (11.04), Ex.\ 11.1]{Olver}.
To prove \eqref{WFAKTA} for $s=\frac 12$ one may use 
\cite[Ch.\ 7 (9.03), (11.02), Ex.\ 11.1]{Olver} to show
$W_{\pm v, c}(Y) = 
O(|c|^{-1} Y^{1/2-|c|})$,
uniformly for $0<Y<\frac 14$ and complex $c$ with
$0 < |c|< \frac 14$ (and fixed $v \in \Z$);
then apply the maximum modulus principle in the $c$-variable along 
$|c| = -(10\log Y)^{-1}$.

Clearly, 
when $s \in \frac 12 + \ii[0, \infty)$ or $s \in (\frac 12,1)$,
\eqref{SMAASSWV10} follows from \eqref{WFAKTA}.

By 
Lemma \ref{HOLO},
the only remaining possibility is $s \in \{1,2,...,|v|\}$,
and in this case we know that $c_n \neq 0$ can only hold if
$nv > 0$. 
Using \eqref{WHITTELEM} and the recursion relation 
$W_{k+1,m}(Y) = (Y/2 -k
- Y (\partial / \partial Y)) W_{k,m}(Y)$
(cf.\ \cite[p.\ 507 (13.4.33)]{AS}),
we now find that:
\begin{equation*}
W_{v \cdot \sgn(n),s - \frac 12}(Y) =
W_{|v|, s - \frac 12} (Y) = O(Y^s)
\qquad \text{as } \: Y \to 0^+.
\end{equation*}
Now \eqref{SMAASSWV10} is completely proved.

Using \eqref{SMAASSWV10} together with the
well-known estimate $\tau(|l|) = O(|l|^{\varepsilon})$,
we obtain
\begin{align} \label{SMAASWV12}
\int_0^1 |S|^2 \, d\alpha \leqq 
|c_0|^2 y^{2-2\sigma} & +
O ( M^{\varepsilon - 1} )
\sum_{n \neq 0}
|c_n|^2 |n|^{\varepsilon-1} e^{-2\pi |n| y} .
\end{align}
Recall here that by 
Lemma \ref{EGENVLEMMA} and Lemma \ref{HOLO},
$c_0 \neq 0$ is possible only if $s \in (\frac 12,1)$.
We have the following mean square bound on the coefficients
$c_n$:
\begin{equation} \label{SMAASSWV11}
\sum_{1 \leqq |n| \leqq N} |c_n|^2 = O(N). 
\end{equation}
(The implied constant depends on $\phi$.)
In the case $\tre s < 1$, \eqref{SMAASSWV11} follows from 
Lemma \ref{MAASSWV} and \cite[Thm.\ 3.2]{IW};
in the remaining case, i.e.\ $s \in \{1,2,...,|v|\}$,
\eqref{SMAASSWV11} follows from
Lemma \ref{HOLO} and \cite[Thm.\ 5.1]{IW2}, and partial summation.

Now \eqref{SMAASSWVRES} follows from \eqref{SMAASWV12} and
\eqref{SMAASSWV11}, by partial summation.
\end{Proof}

\begin{corollary} \label{PHINACOR}
Let $\phi,\lambda,s$ be as in Proposition \ref{SMAASSWV},
and let ${\minexp} > 0$ and $c>0$.
Let $A>1$ be an arbitrary integer.
If $\frac 12 < s < 1$ then we also assume 
$A > {\minexp}^{-1}(2-2s)^{-1}$.
Then there is a set $P=P(\Gamma,\phi,{\minexp},c,A) \subset \R$ 
of full Lebesgue measure such that for each $\alpha \in P$, we have
\begin{equation} \label{PHINACORRES}
\lim_{M \to \infty} \frac 1{M^A} \sum_{m=1}^{M^A} 
\phi \bigl( m \alpha + \ii cM^{-A{\minexp}} \bigr)
= 0. 
\end{equation}
\end{corollary}

\begin{Proof}
Let us fix $\varepsilon > 0$ so small that 
$A (1-({\minexp}+1)\varepsilon) > 1$.
Applying Proposition \ref{SMAASSWV} with $y=cM^{-{\minexp}}$ we obtain
\begin{equation} \label{FSUMBOUND}
\int_0^1 \Bigl |  
\frac 1M \sum_{m=1}^M \phi(m \alpha + \ii cM^{-{\minexp}})
\Bigr |^2 \, d\alpha = 
O (M^{-b}) ,
\qquad
\text{as } \: M \to \infty.
\end{equation}
Here $b = \min(1-({\minexp}+1) \varepsilon, {\minexp}(2-2s))$ 
if $\frac 12<s<1$, otherwise
$b=1-({\minexp}+1) \varepsilon$. It follows from our assumptions that 
$A b > 1$, and hence
\begin{align*}
\int_0^1 \sum_{M=1}^{\infty} \Bigl |  
\frac 1{M^A} \sum_{m=1}^{M^A} 
\phi \bigl( m \alpha + \ii c(M^A)^{-{\minexp}} \bigr)
\Bigr |^2 \, d\alpha 
= O(1) \sum_{M=1}^{\infty} (M^A)^{-b} < \infty.
\end{align*}
It follows that the integrand in the left hand side is finite for
almost all $\alpha$, and hence, a fortiori, there is a set
$P \subset [0,1]$ of full Lebesgue measure such that
\eqref{PHINACORRES} holds for all $\alpha \in P$.

Finally, we notice that the sum 
$\sum_{m} \phi \bigl( m \alpha + \ii y \bigr)$ is invariant under
$\alpha \mapsto \alpha+1$, since $\matr 1101 \in \Gamma$.
The desired result follows from this.
\end{Proof}

\begin{proposition} \label{SEISENST}
Let $\varepsilon>0$, $k \in \{1,...,\kappa\}$, $v \in \Z$,
and take $s = \sigma + it$ either on $\frac 12 + \ii[0,\infty)$
or on $(\frac 12,1)$, such that $E_k(z,s,2v)$ does not have a pole at $s$.
We write
\begin{equation*}
S = S(M,y, \alpha) = \frac 1M \sum_{m=1}^M 
E_k(m \alpha + {\ii y},s,2v) .
\end{equation*}
We then have, for all $M \geqq 1$ and $0 < y \leqq 1$,
\begin{equation*} 
\int_0^1 |S|^2 \, d\alpha = 
O \left( y^{1-2\sigma-\varepsilon} M^{\varepsilon -1} 
+ y^{2-2\sigma} \right) 
\end{equation*}
(The implied constant depends only on $\Gamma, v, s, \varepsilon$.)
\end{proposition}

\begin{Proof}
Mimicing the proof of Proposition \ref{SMAASSWV},
we obtain \eqref{SMAASWV12} with 
the term $|c_0|^2 y^{2-2\sigma}$ replaced by 
$| \delta_{k1} y^s + c_0 y^{1-s} |^2$
--
the $c_n$'s are now the Fourier coefficients in the expansion
of $E_k(z,s,2v)$
as in \eqref{EFOURIEREXP} in Lemma \ref{EISENSTEINRS}.
Using the bound from Lemma \ref{EISENSTEINRS}
in the form $\sum_{1 \leqq |n| \leqq N} |c_n|^2 = O(N^{2 \sigma +
\varepsilon})$,
together with partial summation, 
we obtain $\int_0^1 |S|^2 \, d\alpha \leqq 
O(M^{\varepsilon - 1} y^{1-2\sigma-2\varepsilon} + y^{2-2\sigma}).$ 
This proves the desired bound, with $2\varepsilon$ in place of
$\varepsilon$.
\end{Proof}

\begin{corollary} \label{EISNACOR}
Let $k,v,s=\sigma+\ii t$ be as in Proposition \ref{SEISENST}
(in particular, $\frac 12 \leqq \sigma < 1$),
and let ${\minexp} > 0$ and $c>0$.
We make the assumption that $\sigma < \frac 12 (1+1/{\minexp})$,
and we let $A$ be an arbitrary integer greater than 
\mbox{$\max((1+{\minexp}-2{\minexp}\sigma)^{-1},
(2{\minexp}-2{\minexp}\sigma)^{-1})$.}
Then there is a set $P=P(\Gamma,k,v,s,{\minexp},c,A) \subset \R$ 
of full Lebesgue measure such that for each $\alpha \in P$, we have 
\begin{equation*}
\lim_{M \to \infty} \frac 1{M^A} \sum_{m=1}^{M^A} 
E_k \bigl( m \alpha + \ii cM^{-A{\minexp}},s,2v \bigr)
= 0. 
\end{equation*}
\end{corollary}

\begin{Proof}
This is similar to the proof of Corollary \ref{PHINACOR},
using Proposition \ref{SEISENST} instead of
Proposition \ref{SMAASSWV}.
\end{Proof}

\section{Averages of eigenfunctions: conditional results}
\label{CONDRESSECSTART}

In this section we will show that,
in Corollaries \ref{PHINACOR} and \ref{EISNACOR},
we may under certain conditions
replace the sequence $1^A,2^A,3^A,...$
by the full sequence $1,2,3,...$,
and allow the constant $c$ to vary over any fixed compact interval.
The main reason for this is that 
the hyperbolic distance between $\ii M^{-A \nu}$ and $\ii (M+1)^{-A \nu}$ 
tends to $0$ as $M \to \infty$.
We will start with the simplest case; that of $\phi$ being a cusp form.

Throughout Section \ref{CONDRESSECSTART} and \ref{MAINTHMPRSE}, 
we will let $\mathbf{D}$ denote a fixed
countable dense subset of $\R^+$.

\begin{lemma} \label{USEDISTANCE}
Let $\phi$ be a Maass waveform of even integer weight,
and assume that $\phi$ is a cusp form.
Let ${\minexp}>0$, $\alpha \in \R$ and $A \in \Z^+$.
We assume that, for all $c \in \mathbf{D}$,
\begin{equation} \label{USEDISTANCEASS}
\lim_{M \to \infty} \frac 1{M^A} \sum_{m=1}^{M^A} 
\phi \bigl( m \alpha + \ii cM^{-A{\minexp}} \bigr)
= 0. 
\end{equation}
Then, for any fixed $0<C_1<C_2$, we have
\begin{equation*}
\frac 1{M} \sum_{m=1}^{M} 
\phi \bigl( m \alpha + \ii y \bigr)
\to 0, 
\end{equation*}
uniformly as $M \to \infty$ and 
$C_1M^{-{\minexp}} \leqq y \leqq C_2M^{-{\minexp}}$.
\end{lemma}

\begin{Proof}
Let $0<C_1<C_2$ be given, and fix some $\varepsilon > 0$.
Since $\phi$ is a cusp form, it is bounded and uniformly continuous over all of $\HH$.
We fix positive numbers $B$ and $\delta_0$ such that $|\phi(z)| \leqq B$ holds for all
$z \in \HH$, and such that
$|\phi(z) - \phi(w)| < \varepsilon$ holds whenever $\delta(z,w) < \delta_0$,
where $\delta(\cdot,\cdot)$ denotes hyperbolic distance.

Let us also fix a finite subset $S \subset [C_1,C_2] \cap \mathbf{D}$ such that for 
each $c \in [C_1,C_2]$ there is at least one $c_1 \in S$ with
$|\log(c/c_1)| \leqq \delta_0/2$.
Because of our assumption \eqref{USEDISTANCEASS} and the fact that $S$ is finite,
we know that for all sufficiently large numbers $M_1$,
\begin{equation} \label{USEDISTANCE1}
\Bigl | 
\frac 1{M_1^A}
\sum_{m=1}^{M_1^A} \phi \left( m\alpha + \ii c_1 M_1^{-A{\minexp}} \right) 
\Bigr | \leqq \varepsilon,
\qquad \forall c_1 \in S.
\end{equation}

We now let an arbitrary number $M \in \Z^+$ be given and define
$M_1 \in \{2,3,...\}$ through $(M_1-1)^A \leqq M < M_1^A$.
We make the assumption that $M$ is so large that \eqref{USEDISTANCE1} holds, 
and so that
\begin{equation} \label{USEDISTANCE2}
\frac M{M_1^A} > \max \Bigl( 1 - \frac{\varepsilon}{2B},
\exp \bigl( - \frac{\delta_0}{2{\minexp}} \bigr) \Bigr).
\end{equation}

Given any $c \in [C_1,C_2]$ we 
pick $c_1 \in S$ such that $|\log(c/c_1)| \leqq \delta_0/2$.
We now have
\begin{align} \label{USEDISTANCE3}
& \Bigl | 
\frac 1M  \sum_{m=1}^{M} \phi \left( m\alpha + \ii cM^{-{\minexp}}\right) 
- \frac 1{M_1^A}
\sum_{m=1}^{M_1^A} \phi \left( m\alpha + \ii c_1 M_1^{-A{\minexp}} \right) 
\Bigr | 
\\ \notag
& \leqq  \frac 1M \sum_{m=1}^M \left |
\phi \left( m\alpha + \ii cM^{-{\minexp}}\right) 
- \phi \left( m\alpha + \ii c_1 M_1^{-A{\minexp}} \right)  \right | 
\\ \notag & \qquad 
+ \Bigl( \frac 1M - \frac 1{M_1^A} \Bigr) M B 
+ \frac 1{M_1^A} \bigl( M_1^A - M \bigr) B.
\end{align}
Here the last two terms are equal, and are each less than $\varepsilon/2$,
by \eqref{USEDISTANCE2}.
Furthermore,
\begin{align*}
& \delta \bigl( m \alpha + \ii cM^{-{\minexp}},
m \alpha + \ii c_1M_1^{-A{\minexp}} \bigr)
=
\biggl | \log \Bigl( \frac{cM^{-{\minexp}}}{c_1M_1^{-A{\minexp}}} \Bigr)
\biggr |
\\
& \qquad \leqq | \log (c/c_1) | + {\minexp} | \log (M/M_1^A)| <
\frac{\delta_0}2 + \frac{\delta_0}2 = \delta_0, 
\end{align*}
again by \eqref{USEDISTANCE2}.
It follows that
$\bigl | \phi \left( m\alpha + \ii cM^{-{\minexp}}\right) - 
\phi \left( m\alpha + \ii c_1M_1^{-A{\minexp}}\right) \bigr | < \varepsilon$
holds for each $m$.
Hence the difference in \eqref{USEDISTANCE3} is less than $2 \varepsilon$,
and using \eqref{USEDISTANCE1} we conclude that
\begin{equation*}
\Bigl | 
\frac 1{M}
\sum_{m=1}^{M} \phi \left( m\alpha + \ii c M^{-{\minexp}} \right) 
\Bigr | < 3 \varepsilon.
\end{equation*}
This holds for all sufficiently large $M$,
and all $c \in [C_1,C_2]$.
This concludes the proof.
\end{Proof}

In order to extend the above lemma to the case when $\phi$ is not a cusp form,
or when $\phi$ is replaced by the Eisenstein series, we need a good bound on those 
contributions to our sum which come from points lying far out in the cusps of
$\qspace$.
This is provided by the following lemma.

\begin{lemma} \label{CUSPLEMMA}
Let ${\minexp}>0$, $\alpha \in \R$, $A \in \Z^+$ and $s_0 \in (\frac 12,1)$ be given,
such that $s=s_0$ is not a pole of any of the Eistenstein series
$E_k(z,s,0)$, $k=1,...,\kappa$.
We assume that, for some constant $K_1 > 0$,
\begin{equation} \label{MAINPROPASS1}
\frac 1{M^A} \Bigl | \sum_{m=1}^{M^A} 
E_k\bigl( m \alpha + \ii M^{-A{\minexp}},s_0,0 \bigr) \Bigr | \leqq K_1
, \quad
\forall k \in \{1,...,\kappa\}, \: M \in \Z^+.
\end{equation}
We let $0<C_1<C_2$ be given.
Then there exists a constant $K_2>0$ such that
\begin{equation} \label{CUSPLEMMARES}
\frac 1M \sum_{m=1}^{M} 
\ygg(m\alpha + \ii cM^{-{\minexp}})^{s_0}
\leqq K_2,
\end{equation}
for all $M \in \Z^+$, $c \in [C_1,C_2]$.
\end{lemma}

\begin{Proof}
We let $G(z) = \sum_{k=1}^{\kappa} E_k(z,s_0,0)$;
this is a continuous, real-valued and $\Gamma$-invariant function on
$\HH$.
It follows from \cite[pp.\ 280(Prop.\ 3.7), 297(F)]{TR2}
that $G(z) \sim \ygg(z)^{s_0}$ 
as $z \in \FF$ approaches any of the cusps.
Hence $B = \inf_{z \in \FF} G(z)$ is finite, i.e., $B> - \infty$,
and there is a positive constant $C_3$ such that,
for the function $G_1(z) = G(z) + |B| + 1$:
\begin{equation} \label{CUSPLEMMA10}
G_1(z) \geqq C_3 \ygg(z)^{s_0},
\qquad \forall z \in \FF.
\end{equation}
Then this inequality actually holds for all $z \in \HH$, since both
sides are $\Gamma$-invariant functions.

It follows from \eqref{MAINPROPASS1} that
\begin{equation*}
\frac 1{M^A} \sum_{m=1}^{M^A} 
G_1\bigl( m \alpha + \ii M^{-A{\minexp}} \bigr) 
\leqq |B| + 1 + \kappa K_1, \qquad \forall M \in \Z^+,
\end{equation*}
and hence there is a positive constant $C_4>0$ such that
\begin{equation} \label{M0MBOUND}
\frac 1{M^A} \sum_{m=1}^{M^A} \ygg
\bigl( m \alpha + \ii M^{-A{\minexp}} \bigr)^{s_0} \leqq C_4,
\qquad \forall M \in \Z^+.
\end{equation}
But the invariant height function $\ygg(z)$ satisfies the following
elementary inequality:
\begin{equation} \label{YGGELEM}
\ygg(z_1) \leqq e^{\delta(z_1,z_2)} \ygg(z_2).
\end{equation}
This inequality follows directly from the definition \eqref{YGGDEF}, 
if we write $\delta_0 = \delta(z_1,z_2)$
and notice that 
we have $\delta(N_kW(z_1),N_kW(z_2)) = \delta_0$
for all $k \in \{1,...,\kappa\}$, $W \in \Gamma$,
and hence
\begin{align*}
\tim N_kW(z_1) & \leqq
\sup \bigl\{ \tim z \mid
z \in \HH, \delta(z,N_kW(z_2)) \leqq \delta_0 \bigr\}
\\
& = e^{\delta_0} \tim N_kW(z_2) \leqq
e^{\delta_0} \ygg(z_2).
\end{align*}
Now let $M \in \Z^+$ and let $c$ be any positive real number.
We choose $M_0 \in \{2,3,4,...\}$ such that
$(M_0-1)^A \leqq M < M_0^A$.
We then have
$1 < M_0^A/M \leqq M_0^A/(M_0-1)^A \leqq 2^A$,
and, for any $x \in \R$,
\begin{equation*}
\delta(x+\ii cM^{-{\minexp}}, x+ \ii M_0^{-A{\minexp}})
= \bigl | \log (cM^{-{\minexp}}/M_0^{-A{\minexp}}) \bigr |
\leqq | \log c | + A{\minexp} \log 2.
\end{equation*}
Hence by \eqref{M0MBOUND} and \eqref{YGGELEM}:
\begin{align} \notag
&\frac 1{M} \sum_{m=1}^{M} \ygg
\bigl( m \alpha + \ii cM^{-{\minexp}} \bigr)^{s_0}
\\ \notag
& \leqq
\frac {M_0^A}{M} \cdot 
\frac 1{M_0^A} \sum_{m=1}^{M_0^A} 
\bigl [ e^{|\log c| + A{\minexp} \log 2} \, \ygg 
\bigl( m \alpha + \ii
M_0^{-A{\minexp}} \bigr) \bigr ] ^{s_0} 
< 2^{(1+ s_0{\minexp})A} e^{s_0 |\log c|} C_4.
\end{align}
The lemma follows from this.
\end{Proof}

Using Lemma \ref{CUSPLEMMA}, we are now able to extend 
Lemma \ref{USEDISTANCE} to the case of arbitrary Maass waveforms.

\begin{proposition} \label{MAASSWVLEMMA}
Let ${\minexp}>0$, $\alpha \in \R$, $A \in \Z^+$,
$s_0 \in (\frac 12,1)$,
and assume that \eqref{MAINPROPASS1} in Lemma \ref{CUSPLEMMA} holds;
also let $\phi$ be a nonconstant Maass waveform of even integer weight,
and assume that, for all $c \in \mathbf{D}$,
\begin{equation} \label{USEDISTANCEASSNY}
\lim_{M \to \infty} \frac 1{M^A} \sum_{m=1}^{M^A} 
\phi \bigl( m \alpha + \ii cM^{-A{\minexp}} \bigr)
= 0. 
\end{equation}
We then have, for any fixed constants $0<C_1<C_2$,
\begin{equation}
\frac 1M \sum_{m=1}^{M} \phi \left( m\alpha + {\ii y} \right) 
\to 0,
\end{equation}
uniformly as $M \to \infty$ and 
$C_1 M^{-{\minexp}} \leqq y \leqq C_2 M^{-{\minexp}}$.
\end{proposition}

\begin{Proof}
If $\phi$ is a cusp form then this was proved in Lemma \ref{USEDISTANCE},
and the assumption \eqref{MAINPROPASS1} is not needed.
We now assume that $\phi$ is not a cusp form.
We write $\phi \in L_2(\lambda,2v),$ $\lambda=s(1-s)$ as in Lemma \ref{EGENVLEMMA}.
We then have $\frac 12 < s < 1$ by Lemma \ref{EGENVLEMMA} and Lemma \ref{HOLO},
and by \eqref{MWGROWTH} there is a constant $C>0$
such that $|\phi(z)| \leqq C \ygg(z)^{1-s}$
for all $z \in \HH$.

We let $0<C_1<C_2$ be given, and fix some $\varepsilon > 0$.
We take $K_2>0$ so that \eqref{CUSPLEMMARES} in Lemma \ref{CUSPLEMMA}
holds, and we fix a constant $Y>0$ so large that
\begin{equation} \label{MAASSWVLEMMANY10}
K_2 C Y^{1-s-s_0} < \varepsilon.
\end{equation}
We let $H: \R^+ \to [0,1]$ be a continuous function such that
$H(y)=1$ for $0<y \leqq Y$ and
$H(y)=0$ for $Y+1 \leqq y$, and then let
\begin{equation} \label{MAASSWVLEMMANY3}
f(z) = H(\ygg(z)) \cdot \phi(z).
\end{equation}
Now, if $\ygg(z) > Y$,
we have $|\phi(z) - f(z)|
\leqq |\phi(z)| \leqq C \ygg(z)^{1-s} \leqq C
Y^{1-s-s_0} \ygg(z)^{s_0}.$ If $\ygg(z) \leqq Y$, the same difference is 
$|\phi(z) - f(z)|=0$.
Hence:
\begin{align*}
|\phi(z) - f(z)| \leqq C Y^{1-s-s_0} \ygg(z)^{s_0}
\qquad \text{for \textit{all} $z \in \HH$.}
\end{align*}
Using this together with 
\eqref{CUSPLEMMARES}
and \eqref{MAASSWVLEMMANY10},
we obtain
\begin{equation} \label{MAASSWVLEMMANY1}
\Bigl | \frac 1M \sum_{m=1}^M \phi(m \alpha + \ii cM^{-{\minexp}}) 
- \frac 1M \sum_{m=1}^M f(m \alpha + \ii cM^{-{\minexp}}) \Bigr | 
< \varepsilon
\end{equation}
for all $M \in \Z^+$, $c \in [C_1,C_2]$.


In particular, for each $c \in \mathbf{D}$,
\eqref{USEDISTANCEASSNY} and \eqref{MAASSWVLEMMANY1} imply
\begin{equation} \label{MAASSWVLEMMANY2}
\limsup_{M \to \infty}
\Bigl | \frac 1{M^A} \sum_{m=1}^{M^A} 
f \bigl( m \alpha + \ii cM^{-A{\minexp}} \bigr) \Bigr | \leqq \varepsilon.
\end{equation}
But it follows from our definition in \eqref{MAASSWVLEMMANY3} 
that $f$ is a function of weight $2v$ (cf.\ \eqref{WEIGHTK})
which has compact support
on $\qspace$; in particular, $f$ is bounded and uniformly continuous on all of $\HH$.
Hence, by arguing as in the proof of Lemma \ref{USEDISTANCE},
using 
\eqref{MAASSWVLEMMANY2} in place of
\eqref{USEDISTANCEASS},
we find that for all sufficiently large $M$,
\begin{equation*}
\Bigl | \frac 1{M} \sum_{m=1}^{M} 
f \bigl( m \alpha + \ii y \bigr) \Bigr | < 4\varepsilon,
\qquad  \forall y \in [C_1 M^{-{\minexp}},C_2 M^{-{\minexp}}].
\end{equation*}
Hence for these $M,y$ we also have, by \eqref{MAASSWVLEMMANY1},
\begin{equation*}
\Bigl | \frac 1{M} \sum_{m=1}^{M} 
\phi \bigl( m \alpha + \ii y \bigr) \Bigr | < 5\varepsilon.
\end{equation*}
This concludes the proof.
\end{Proof}

\begin{proposition} \label{EAPRLEMMA}
Let ${\minexp}>0$, $\alpha \in \R$, $A \in \Z^+$, $s_0 \in (\frac 12,1)$,
and assume that \eqref{MAINPROPASS1} in Lemma \ref{CUSPLEMMA} holds;
also let $v \in \Z$, $k \in \{1,...,\kappa \}$,
and assume that
\begin{equation} \label{MAINPROPASS2}
\lim_{M \to \infty} \frac 1{M^A} \sum_{m=1}^{M^A} 
E_k \bigl( m \alpha + \ii cM^{-A{\minexp}}, \sfrac 12+it,2v \bigr)
= 0 
\end{equation}
holds for all $c,t \in \mathbf{D}$.
We then have, for any fixed constants $0<C_1<C_2$ and $T>0$,
\begin{equation}
\frac 1M \sum_{m=1}^{M} 
E_k \bigl( m\alpha + {\ii y},\sfrac 12 + \ii t, 2v \bigr) 
\to 0,
\end{equation}
uniformly as $M \to \infty$, $C_1 M^{-{\minexp}} 
\leqq y \leqq C_2 M^{-{\minexp}}$
and $0 \leqq t \leqq T$.
\end{proposition}

\begin{Proof}
We first note that
by \cite[pp.\ 280 (Prop.\ 3.7), 297(F), 368 (5.21),
374 (item 15)]{TR2},
there is a constant $C_5= C_5(\Gamma,T,v)>0$ such that
\begin{equation} \label{PHIKONV1}
\bigl | E_k(z, \sfrac 12+ \ii t,2v) \bigr |
\leqq C_5 \ygg(z)^{{\frac 12}},
\qquad \forall t \in [0,T], \: z \in \HH.
\end{equation}
Hence, since $\frac 12-s_0 < 0$,
we may imitate the proof of
Proposition \ref{MAASSWVLEMMA}, using \eqref{MAINPROPASS2}, to prove that
for any \textit{fixed $t \in \mathbf{D}$,}
\begin{equation} \label{PHIKONV3}
\frac 1M \sum_{m=1}^{M} 
E_k \bigl( m\alpha + {\ii y},\sfrac 12 + \ii t, 2v \bigr) 
\to 0,
\end{equation}
uniformly as $M \to \infty$ and 
$C_1 M^{-{\minexp}} \leqq y \leqq C_2 M^{-{\minexp}}$.

Now let $\varepsilon > 0$ be given.
We let $K_2$ be as in Lemma \ref{CUSPLEMMA}.
We claim that there is a number $t_0 > 0$ such that
for all $z \in \HH$ and all $t_1 , t_2 \in [0,T]$,
\begin{equation} \label{PHIKONV2}
|t_1-t_2| < t_0 \Longrightarrow
\Bigl | E_k(z,\sfrac 12 + \ii t_1, 2v) - 
E_k(z,\sfrac 12 + \ii t_2, 2v) \Bigr | \leqq
\frac{\varepsilon}{K_2} \ygg(z)^{s_0}
\end{equation}
To prove this, first note that we may assume $z \in \FF$,
by $\Gamma$-invariance.
We take $Y > 0$ so large that
$2C_5 Y^{\frac 12 - s_0} \leqq \varepsilon / K_2$;
then \eqref{PHIKONV2} holds automatically for all
$z \in \FF$ with $\ygg(z) \geqq Y$, by \eqref{PHIKONV1}.
But the region 
$\FF \cap \{ \ygg(z) \leqq Y \}$ is \textit{compact} 
(as is the interval $[0,T]$),
and $E_k(z,\sfrac 12+ \ii t, 2v)$ is a continuous function of 
$\langle z,t \rangle.$
Hence we can indeed choose $t_0 > 0$ so small that
\eqref{PHIKONV2} holds.

Next, we fix a finite subset $S \subset [0,T] \cap \mathbf{D}$
such that for each $t \in [0,T]$ there is at least one 
$t_1 \in S$ with $|t-t_1| < t_0$.
It follows from \eqref{PHIKONV3} that there is a number $M_0$ such that, 
for all integers $M \geqq M_0$
and all $y \in [C_1 M^{-{\minexp}},C_2 M^{-{\minexp}}],$ $t_1 \in S$,
\begin{equation} \label{PHIKONV10} 
\Bigl | \frac 1M \sum_{m=1}^{M} 
E_k \bigl( m\alpha + {\ii y},\sfrac 12 + \ii t_1, 2v \bigr)  \Bigr | 
\leqq \varepsilon.
\end{equation}

Now, for $t \in [0,T]$ arbitrary, we take $t_1 \in S$ such that
$|t-t_1| < t_0$.
It then follows from \eqref{PHIKONV10}, \eqref{PHIKONV2} 
and 
\eqref{CUSPLEMMARES} that, for all $M \geqq M_0$
and all $y \in [C_1 M^{-{\minexp}},C_2 M^{-{\minexp}}]$,
\begin{equation*}
\Bigl | \frac 1M \sum_{m=1}^{M} 
E_k \bigl( m\alpha + {\ii y},\sfrac 12 + \ii t, 2v \bigr)  \Bigr | 
\leqq 2\varepsilon.
\end{equation*}
This concludes the proof.
\end{Proof}

\section{Equidistribution for unbounded test functions}
\label{MAINTHMPRSE}

In this section we will prove Theorem \ref{MAINTHM},
using the results from the preceding section
together with the spectral expansion 
\eqref{THFOURIER}, \eqref{SPEKTRALUTV}.
In fact, we will prove the slightly stronger Theorem \ref{MAINTHM}' below,
which extends \mbox{Theorem \ref{MAINTHM}}
to cases of unbounded test functions.

\label{PSEDEFPAGE}
$\quad$\textbf{Definition.} 
Given $\gamma \geqq 0$, we let $B_{\gamma}(\Tqspace)$ be the family of 
continuous $\Gamma$-invariant functions $f:T_1 \HH \to \C$
such that $| f(z,\theta) | \leqq C \ygg(z)^{\gamma}$
holds for all $(z,\theta) \in T_1\HH$ and some constant $C>0$.
Given ${\minexp}>0$ and $\alpha \in \R$,
we say that $\langle \Gamma,{\minexp},\alpha \rangle$-PSE$_{\gamma}$ holds
($\langle \Gamma,{\minexp},\alpha \rangle$-Point Set Equidistribution)
if, for any fixed $f \in B_{\gamma}(\Tqspace)$ 
and any fixed numbers $0<C_1<C_2$,
we have
\begin{equation} \label{MAINTHMRESDEF}
\frac 1M \sum_{m=1}^{M} f \left( m\alpha + {\ii y} ,0 \right) 
\to \langle f \rangle,
\end{equation}
uniformly as $M \to \infty$ and
$C_1M^{-{\minexp}}\leqq y \leqq C_2M^{-{\minexp}}$.
\vspace{5pt}

Hence, Theorem \ref{MAINTHM} says that for any given ${\minexp}>0$,
$\langle \Gamma,{\minexp},\alpha \rangle$-PSE$_0$ holds for almost 
all $\alpha \in \R$.
\vspace{5pt}

\enlargethispage{25pt}
$\: \:$\textbf{Theorem \ref{MAINTHM}'.}
\textit{Let ${\minexp}>0$.
Then there is a set $P=P(\Gamma,{\minexp}) \subset \R$ of 
full Lebesgue measure such
that for each $\alpha \in P$ and each
\begin{equation} \label{MAINTHMPASS}
0 \leqq \gamma < \min(1,\sfrac 12(1+1/{\minexp})),
\end{equation}
$\langle \Gamma,{\minexp},\alpha \rangle$-PSE$_{\gamma}$ holds.}

\begin{remark}
It is essential that the summation in \eqref{MAINTHMRESDEF} does not
include $m=0$.
Indeed, if $\Gamma$ has a cusp at the point $0$,
then for $f \in B_{\gamma}(\Tqspace)$, we may have
$f({\ii y},0)$ increasing like $y^{-\gamma}$ as $y \to 0^+$,
thus causing $\frac 1M f(iy,0) \to \infty$ as $M \to \infty$,
if ${\minexp} \gamma > 1$.
An example where ${\minexp} \gamma = 1$ and where the contribution for $m=0$
has to be treated separately is found in the proof of 
\mbox{Proposition \ref{APPLPROP}} below.

On the other hand, if ${\minexp} \gamma < 1$ then for any $\Gamma$ and 
any fixed $x \in \R$
we have $\frac 1M f(x+{\ii y},0) \to 0$ as $M \to \infty$ and
$C_1 M^{-{\minexp}} \leqq y \leqq C_2 M^{-{\minexp}}$.
(To see this, write $z=x+iy$ and $\ygg(z) = \tim N_k W(z)$ with 
$W \in \Gamma$; for $y$ small one then knows that $|c| \geqq 1$ in
$\matr abcd = N_k W$, and this forces
$|f(z,0)| \leqq C \ygg(z)^{\gamma}
\leqq C(c^{-2} y^{-1})^{\gamma} \leqq C y^{-\gamma}.$)
Hence when $\nu \gamma < 1$, \mbox{Theorem \ref{MAINTHM}'} remains true if
$\sum_{m=1}^M$  in \eqref{MAINTHMRESDEF}
is replaced by (say) $\sum_{m=0}^M$ or $\sum_{m=0}^{M-1}$.
\end{remark}

\label{MAINPROPPROOF3}
\textbf{Proof of Theorem \ref{MAINTHM}'.} 
Since any countable intersection of sets of full Lebesgue measure
is again of full Lebesgue measure, it suffices to prove that
for any \textit{fixed} $\gamma$ satisfying \eqref{MAINTHMPASS},
there is a subset $P = P(\Gamma,{\minexp},\gamma) \subset \R$ such that
$\langle \Gamma,{\minexp},\alpha \rangle$-PSE$_{\gamma}$ 
holds for all $\alpha \in P$.

From now on, we keep ${\minexp},\gamma$ fixed as above.
We choose a number $s_0$ with
$\max(\frac 12,\gamma) < s_0 < \min(1,\frac 12(1+1/{\minexp}))$ such that
$s=s_0$ is not a pole of any of the Eistenstein series
$E_k(z,s,0)$, $k=1,...,\kappa$.
For each $v \in \Z$, we let $\mathbf{D}_{2v}$ be a complete set of discrete
eigenfunctions of $- \Delta_{2v}$, taken to be orthonormal.
Using Corollaries \ref{PHINACOR} and \ref{EISNACOR},
and the fact that any countable intersection of sets of full Lebesgue measure
is again of full Lebesgue measure, it follows that there 
exists a positive integer $A=A(\Gamma,{\minexp},s_0)$ 
and a set $P \subset \R$ of full Lebesgue measure 
such that for all $\alpha \in P$,
the assumption in Lemma \ref{CUSPLEMMA} is fulfilled,
as well as the assumption in Proposition \ref{MAASSWVLEMMA}
for each nonconstant $\phi \in \cup_v \mathbf{D}_{2v}$,
and the assumption in Proposition \ref{EAPRLEMMA} for all $v \in \Z$, $k \in \{1,...,\kappa\}$.

\label{MAINPROPPROOF4}
Now take any $\alpha \in P$, and fix some numbers $0<C_1<C_2$.
We are going to prove that for each fixed $f \in B_{\gamma}(\Tqspace)$,
we have
\begin{equation} \label{MAINTHVILL}
\frac 1M \sum_{m=1}^{M} f \left( m\alpha + {\ii y} ,0 \right) 
\to \langle f \rangle,
\end{equation}
uniformly as $M \to \infty$ and
$C_1M^{-{\minexp}}\leqq y \leqq C_2M^{-{\minexp}}$.
This will complete the proof of Theorem \ref{MAINTHM}'.

Let us define a norm on $B_{\gamma}(\Tqspace)$ through
\begin{equation*}
|| f ||_{(s_0)}
:= \sup_{(z,\theta) \in T_1\HH} 
\, \Bigl( |f(z,\theta)| \, \ygg(z)^{-s_0} \Bigr).
\end{equation*}
It follows from the definition of $B_{\gamma}(\Tqspace)$
that $|| f ||_{(s_0)}$ is finite for each 
$f \in B_{\gamma}(\Tqspace)$,
since $\gamma < s_0$ and $\ygg(z)$ is bounded from below by a positive constant.

Now by Lemma \ref{CUSPLEMMA}, we have
\begin{align*}
\Bigl | \frac 1M \sum_{m=1}^{M} f_1 \left( m\alpha + {\ii y} ,0 \right) 
-
\frac 1M \sum_{m=1}^{M} f_2 \left( m\alpha + {\ii y} ,0 \right) 
\Bigr | \leqq K_2 || f_1-f_2 ||_{(s_0)},
\end{align*}
for all $f_1,f_2 \in B_{\gamma}(\Tqspace)$
and all $M \in \Z^+$, $y \in [C_1M^{-{\minexp}}, C_2M^{-{\minexp}}]$.
Furthermore, there is a constant $C_6>0$ which only depends on
$\Gamma$ and $s_0$ such that
\begin{equation*}
\Bigl | \langle f_1 \rangle - \langle f_2 \rangle \Bigr |
\leqq C_6 || f_1-f_2 ||_{(s_0)},
\qquad
\forall f_1,f_2 \in B_{\gamma}(\Tqspace)
\end{equation*}
(for recall  \eqref{YGGINTEGRAL} and $s_0<1$).
Because of these last two inequalities,
it now suffices to prove \eqref{MAINTHVILL}
for all $f \in S$,
where $S$ is any \textit{$|| \cdot ||_{(s_0)}$-dense subset}
of $B_{\gamma}(\Tqspace)$.

We claim that $C_c(\Tqspace)$
is $|| \cdot ||_{(s_0)}$-dense in $B_{\gamma}(\Tqspace)$.
For take $f \in B_{\gamma}(\Tqspace)$ arbitrary.
Given $Y \geqq 100$ we define $H:\R^+ \to [0,1]$ as in the proof
of Proposition \ref{MAASSWVLEMMA}, and let
$f_1(z,\theta) = H(\ygg(z)) \cdot f(z,\theta).$
Then $f_1 \in C_c(\Tqspace)$.
Arguing as in the proof of Proposition \ref{MAASSWVLEMMA},
we find that $|| f - f_1 ||_{(s_0)}$ can be made
arbitrarily small by taking $Y$ sufficiently large.
This proves our claim.

Next, by a standard convolution argument,
the space \mbox{$C^{\infty}(T_1\HH) \cap C_c(\Tqspace)$}
is dense in $C_c(\Tqspace)$ with respect to the usual 
supremum norm, and hence, a fortiori, with respect to
$|| \cdot ||_{(s_0)}.$

Now let $f \in C^{\infty}(T_1\HH) \cap C_c(\Tqspace)$ be given.
We expand $f$ in a Fourier series as in \eqref{THFOURIER}.
Integrating by parts twice in \eqref{THFOURIERCOEFF} we find that
$|\widehat{f}_v(z)| \leqq (2\pi)^{-1} |v|^{-2} \int_0^{2\pi}
\bigl | (\partial^2 / \partial \theta^2) f(z,\theta) \bigr | \, d\theta$
for all $|v| \geqq 1$.
Hence, since $f$ is smooth and of compact support modulo $\Gamma$,
we have uniform convergence in \eqref{THFOURIER},
over all of $T_1 \HH$.
In particular, the space of \textit{finite} sums
$\sum_v f_v(z) e^{iv \theta}$, with 
$f_v \in C^{\infty}(\HH) \cap C_c(\qspace,2v)$,
is $|| \cdot ||_{(s_0)}$-dense in $B_{\gamma}(\Tqspace)$.

Hence, by linearity, it will now be sufficient to prove 
\eqref{MAINTHVILL} for functions $f$ of the form
$f(z,\theta) = f_v(z) e^{iv \theta}$,
where $v \in \Z$ and $f_v \in C^{\infty}(\HH) \cap C_c(\qspace,2v)$.
For such $f$, we clearly have $\langle f \rangle = 0$
if $v \neq 0$ (since $\int_0^{2\pi} f(z,\theta) \, d\theta = 0$ for all
$z$), and $\langle f \rangle = \mu(\FF)^{-1} \int_{\FF}
f_v(z) \, d\mu(z)$ if $v=0$.
Hence, our goal is now to prove, for any fixed $v \in \Z$, and any
fixed $f_v \in C^{\infty}(\HH) \cap C_c(\qspace,2v)$:
\begin{equation} \label{MAINTHVILL2}
\frac 1M \sum_{m=1}^{M} f_v ( m\alpha + {\ii y} ) 
\to \delta_{v0} \, \mu(\FF)^{-1} \int_{\FF} f_v(z) \, d\mu(z),
\end{equation}
uniformly as $M \to \infty$ and 
$C_1M^{-{\minexp}}\leqq y \leqq C_2M^{-{\minexp}}$.

\enlargethispage{10pt}

We will now apply the spectral expansion \eqref{SPEKTRALUTV}.
By an approximation argument of the same type as above,
using the analog of the $|| \cdot ||_{(s_0)}$-norm
for functions $\HH \to \C$ of weight $2v$,
and using Lemma \ref{HUGGSPEKTRAL}, the inequality 
$1+\sqrt{\ygg(z)} \leqq O(\ygg(z)^{s_0})$
and Lemma \ref{CUSPLEMMA},
we find that it suffices to prove
\eqref{MAINTHVILL2} for $f_v$
in the family of \textit{finite} spectral expansions;
\begin{equation*}
f_v(z) = \sum_{\lambda_n \leqq 1/4 + T^2} d_n \phi_n(z)
+  \sum_{k=1}^{\kappa} \int_0^T g_k(t) E_k(z,\sfrac 12 + it,2v) \, dt,
\end{equation*}
where $\phi_n \in {\mathbf D}_{2v}$,
$T \geqq 1$, $d_0,d_1,...$ are arbitrary complex numbers,
and $g_k(t)$ are arbitrary continuous functions on $[0,T]$.
But for \textit{such} a function $f_v$,
Proposition \ref{EAPRLEMMA} and Proposition \ref{MAASSWVLEMMA} can be applied
directly.
Concerning  Proposition \ref{MAASSWVLEMMA}, we
notice that if $v \neq 0$ then all $\phi_n$'s are nonconstant;
if $v=0$ then $\phi_n$ is nonconstant for all $n \neq 0$,
while $\phi_0 \equiv \mu(\FF)^{-\frac 12}$.
Hence we obtain:
\begin{equation*}
\frac 1M \sum_{m=1}^{M} f_v ( m\alpha + {\ii y} ) 
\to \delta_{v0} \, d_0 \mu(\FF)^{-\frac 12},
\end{equation*}
uniformly as $M \to \infty$ and 
$C_1M^{-{\minexp}}\leqq y \leqq C_2M^{-{\minexp}}$.
But here, if $v =0$, we have $d_0 = \langle f_0, \phi_0 \rangle$,
and hence $ d_0 \mu(\FF)^{-\frac 12}=
\mu(\FF)^{-1} \int_{\FF} f_0(z) \, d\mu(z)$.
Hence \eqref{MAINTHVILL2} holds.
This concludes the proof of Theorem \ref{MAINTHM}'.
$\square$ $\square$ $\square$
\label{MAINPROPPROOF2}

\section{Proof of Theorem \ref{MAINTHM2}}

For $x \in \R$ we use $|| x ||$ to denote the distance from $x$ to the
closest integer, i.e. $|| x || = \inf_{n \in \Z} |x-n|.$

\begin{lemma} \label{DIOFSUMLEMMA}
Let $\alpha \in \R$ be of type $K \geqq 2$.
We then have, for all integers $M \geqq 1,$
$N_2 \geqq N_1 \geqq 1$:
\begin{equation} \label{DIOFSUMLEMMARES}
\sum_{n=N_1}^{N_2} \min \Bigl( M, \frac 1{|| n \alpha ||} \Bigr)
\leqq O \Bigl( N_2^{K-1} \log (2N_2) \Bigr).
\end{equation}
Writing $N=N_2-N_1$, the same sum is also
\begin{equation} \label{DIOFSUMLEMMARES2}
\qquad \leqq O \Bigl(
M + \bigl( ( NM ) ^{\frac{K-1}K} + N \bigr) \log (NM+1) \Bigr).
\end{equation}
The implied constants
depend only on $\alpha,K$, not on $M,N_1,N_2$.
\end{lemma}

We remark that whenever $N_2 \geqq M^{\frac 1{K-1}}$ holds,
the bound in \eqref{DIOFSUMLEMMARES2}
is better than (or at least as good as)
the bound in \eqref{DIOFSUMLEMMARES}.

\begin{Proof}
By assumption, there exists a constant $C \in (0,1)$ such that
$| \alpha - a/q | > Cq^{-K}$ 
for all $a \in \Z$, $q \in \Z^+$. 

We first prove \eqref{DIOFSUMLEMMARES2}. 
The bound is trivial when $N=0$,
so we may now assume $N \geqq 1$. 
By Dirichlet's Theorem, given any real number $Q \geqq 1$
there exist integers $a,q$ with $1 \leqq q \leqq Q$,
$\text{gcd}(a,q)=1$ and $| \alpha - a/q | < (qQ)^{-1}
\leqq q^{-2}$.
Using $| \alpha - a/q | > Cq^{-K}$ we then have
$(CQ)^{\frac 1{K-1}} < q \leqq Q$.

Splitting the summation range $N_1 \leqq n \leqq N_2$ into
consecutive blocks of the form
$\{hq+r \mid 1 \leqq r \leqq q \}$
(possibly overshooting on both ends),
and using \cite[Lemma 4.9]{Nathanson} for each block, we now have
\begin{align*}
& \sum_{n=N_1}^{N_2} \min \Bigl( M, \frac 1{|| n \alpha ||} \Bigr)
\leqq O \Bigl( \Bigl( \frac{N_2-N_1}q + 2 \Bigr) (M + q \log q) \Bigr)
\\
& = O \Bigl( \frac{NM}q + M + (N+q) \log q \Bigr)
= O \Bigl( NM Q^{-\frac 1{K-1}} + M + (N+Q) \log Q \Bigr).
\end{align*}
Taking $Q = (NM)^{\frac {K-1}K}$
(remembering $N \geqq 1$),
we obtain \eqref{DIOFSUMLEMMARES2}.

To prove \eqref{DIOFSUMLEMMARES},
we instead let
$Q=C^{-1}(2N_2)^{K-1} > 1$, and find integers
$a,q$ as above. 
Now $Q \geqq q > (CQ)^{\frac 1{K-1}} = 2N_2$,
so by \cite[Lemma 4.8]{Nathanson},
\begin{align*}
\sum_{n=N_1}^{N_2} \min \Bigl( M, \frac 1{|| n \alpha ||} \Bigr)
\leqq \sum_{1 \leqq n \leqq q/2} \frac 1{|| n \alpha ||}
= O \bigl( q \log q \bigr) = O \bigl( N_2^{K-1} \log (2N_2) \bigr).
\end{align*}
\end{Proof}

\begin{proposition} \label{DIOFLEMMA1}
Let $\beta,\alpha,K,{\minexp}$
be as in Theorem \ref{MAINTHM2}.
Let $\phi$ be a non-constant Maass waveform of even integer
weight.
We then have
\begin{equation} \label{DIOFLEMMA1ASS}
\frac 1M \sum_{m=1}^M \phi(m \alpha + {\ii y})
\to 0,
\end{equation}
uniformly as $M \to \infty$, $y \to 0^+$ 
so long as $y \geqq M^{-{\minexp}}$.
\end{proposition}

\begin{Proof}
We write $\phi \in L_2(\lambda,2v)$ and $\lambda = s(1-s)$
as in Lemma \ref{EGENVLEMMA}, and $s = \sigma + it$.
As in the proof of \mbox{Proposition \ref{SMAASSWV}},
the sum in \eqref{DIOFLEMMA1ASS} is equal to
\begin{equation} \label{DIOFLEMMA10}
c_0 y^{1-s} +
\frac 1M \sum_{n \neq 0} 
\frac {c_n}{\sqrt{|n|}} W _{v \cdot \sgn(n) ,s - \frac 12} 
(4 \pi |n| y) \sum_{m=1}^M e(nm\alpha).
\end{equation}
We will always keep $y<1$.
Recall that we have either $\sigma < 1$,
or else $s=\sigma \geqq 1$ and $c_0 = 0$ (cf.\ Lemma \ref{HOLO}).
Hence we certainly have $c_0 y^{1-s} \to 0$ as $y \to 0^+$.
We now turn to the sum over $n\neq 0$ in \eqref{DIOFLEMMA10}.
We will first work under the assumption that either
$\sigma \leqq \frac 12+\beta$ or $s=\sigma\geqq 1$. 
By reviewing the proof of \eqref{SMAASSWV10} on p.\
\pageref{SMAASSWV10},
we then find that we have the following convenient bound
for all $Y>0$, in all non-vanishing 
terms in \eqref{DIOFLEMMA10}:
\begin{equation} \label{DIOFLEMMA20}
W_{v \cdot \sgn(n) ,s - \frac 12} (Y) =
O(Y^{\frac 12 - \beta} e^{-Y/4}).
\end{equation}
Furthermore, it follows from the definition of $\beta$
and Lemma \ref{MAASSWV} or \mbox{Lemma \ref{HOLO}} that
\begin{equation*}
|c_n| = O(|n|^{\beta+\varepsilon}), \qquad
\text{for all $n \neq 0$.}
\end{equation*}
(Here and in all later ``big-$O$'' estimates in this proof,
the implied constant depends on $\phi$ and $\varepsilon$.)
We also note that
\begin{equation*}
\Bigl | \sum_{m=1}^M e(nm\alpha) \Bigr | = 
\biggl | \frac{e(Mn \alpha) - 1}
{e(n\alpha)-1} \biggr |
\leqq \min(M, || n \alpha ||^{-1}),
\end{equation*}
since $|e(x)-1| = 2|\sin(\pi x)| \geqq 2 ||x||$ for all $x \in \R$.
Hence \eqref{DIOFLEMMA10}
minus the trivial term $c_0 y^{1-s}$ is
\begin{align} \notag
& = 
O \bigl( M^{-1} y^{\frac 12 - \beta} \bigr) \sum_{n=1}^{\infty}
e^{-\pi n y} n^{\varepsilon}  
\min \bigl( M, || n \alpha ||^{-1} \bigr)
\\ \label{DIOFLEMMA11}
& =
O \bigl( M^{-1} y^{\frac 12 - \beta} \bigr) 
\sum_{k=0}^{\infty}
e^{-\pi k} \Bigl( \frac{k+1}y \Bigr)^{\varepsilon}
\sum_{k/y < n \leqq (k+1)/y}
\min \bigl( M, || n \alpha ||^{-1} \bigr).
\end{align}
Applying the bound
\eqref{DIOFSUMLEMMARES2} in Lemma \ref{DIOFSUMLEMMA},
and using $y \geqq M^{-{\minexp}}$, we get:
\begin{align} \notag
& 
= 
O \bigl( M^{-1} y^{\frac 12 - \beta - \varepsilon} \bigr) 
\sum_{k=0}^{\infty}
e^{-3 k}
\bigl( M + (M/y)^{\frac{K-1}K+\varepsilon} 
+ M^{\varepsilon}y^{-1-\varepsilon} \bigr)
\\ \label{DIOFLEMMA12}
&
= O \bigl( y^{\frac 12 - \beta - \varepsilon} 
+ M^{-\frac 1K + \varepsilon + 
       \nu(\frac 12 + \beta - \frac 1K + 2\varepsilon)} 
+ M^{\varepsilon-1 + \nu (\frac 12 + \beta + 2\varepsilon)}
\bigr).
\end{align}
Clearly, if $\beta < \frac 12$, it is possible to keep $\varepsilon$
so small that $\frac 12 - \beta - \varepsilon > 0$;
then the first term above tends to $0$ as $M \to \infty$, $y \to 0^+$.
One also checks that when $\beta < \frac 12$,
\eqref{MAINTHM2ASS} implies that \textit{both}
$\nu < \frac 2{1+2\beta}$
and $\nu < \frac 2{2K\beta + K - 2}$ hold;
hence the last two terms above tend to $0$ as well,
provided that $\varepsilon$ is sufficiently small.

On the other hand, if $\beta \geqq \frac 12$,
we instead apply the bound
\eqref{DIOFSUMLEMMARES} 
(in the crude form ``$O(N_2^{K-1+\varepsilon})$'');
this gives that \eqref{DIOFLEMMA11} is
\begin{align*}
& 
= 
O \bigl( M^{-1} y^{\frac 12 - \beta - \varepsilon} \bigr) 
\sum_{k=0}^{\infty} e^{-3 k} \Bigl( \frac {k+1}y \Bigr)^{K-1+\varepsilon}
= 
O \bigl( M^{-1 + {\minexp}(K + \beta - \frac 32 + 2\varepsilon)} \bigr) .
\end{align*}
By \eqref{MAINTHM2ASS},
this tends to $0$ as $M \to \infty$, provided that $\varepsilon$ is 
sufficiently small.

We now turn to the remaining case; $\frac 12 + \beta < \sigma < 1$.
(Notice that $0 < \beta < \frac 12$ must hold in this case.)
Let us write
\begin{equation}
B(X) = \sum_{1 \leqq n \leqq X} \min \bigl(  M, || n \alpha ||^{-1} \bigr).
\end{equation}
The bound \eqref{DIOFLEMMA20} still holds as $Y \to \infty$, and
thus all terms with $|n| > 1/y$ in \eqref{DIOFLEMMA10} can
be treated exactly as in \eqref{DIOFLEMMA11}, \eqref{DIOFLEMMA12}.
When $Y \to 0^+$, 
we have $W_{v \cdot \sgn(n) ,s - \frac 12} (Y) = O(Y^{1-\sigma})$
(cf.\ \eqref{WFAKTA}).
Hence we obtain the following bound on the 
remaining part (``$0<|n| \leqq 1/y$'') of \eqref{DIOFLEMMA10}:
\begin{align*}
& O \bigl( M^{-1} y^{1- \sigma} \bigr) 
\int_{1/2}^{1/y}
X^{\beta + \frac 12 - \sigma + \varepsilon} \, dB(X)
\\
& = O \bigl( M^{-1} y^{1- \sigma} \bigr) \left \{
\Bigl [ X^{\beta+\frac 12 - \sigma+ \varepsilon} 
B(X) \Bigr ] ^{X=1/y}_{X=1/2}
+  \int^{1/y}_{1/2}
X^{\beta -\frac 12 - \sigma + \varepsilon} B(X) \, dX \right \} .
\end{align*}
Here $B(X) = 0$ for $0<X < 1$, and
it follows from Lemma \ref{DIOFSUMLEMMA} that
$B(X) = O \bigl( X^{K-1+\varepsilon} \bigr)$ for all $X \geqq 1$,
and also, by \eqref{DIOFSUMLEMMARES2},
$B(X) = O \bigl(
(XM)^{\frac{K-1}K+\varepsilon} + X^{1+\varepsilon} \bigr)$ 
whenever $X > M^{\frac 1{K-1}}$; hence we get:
\begin{align*}
\leqq &  O \bigl( M^{-1} y^{1- \sigma} \bigr) \biggl \{
\bigl( M^{\frac 1{K-1}} \bigr)^{\beta 
+ \frac 12 - \sigma + K-1+ 2\varepsilon}
+ y^{-\beta-\frac 12 + \sigma - \varepsilon} 
\bigl((M/y)^{\frac{K-1}K+\varepsilon} + y^{-1-\varepsilon}  \bigr)
\\
& 
+ \int_1^{M^{\frac 1{K-1}}}  X^{\beta - \frac 32 - \sigma + K+2\varepsilon}
\, dX
+ \int_1^{1/y} X^{\beta-\frac 12 - \sigma+\varepsilon}
\bigl((XM)^{\frac{K-1}K + \varepsilon} + X^{1+\varepsilon} \bigr) \, dX
\biggr \}.
\end{align*}
Here all the total $X$-exponents in the integrals are $> -1$,
and hence we obtain:
\begin{align*}
& = O \biggl(
M^{\frac{\beta+1/2-\sigma+2\varepsilon}{K-1}} y^{1-\sigma}
+
M^{-\frac 1K+\varepsilon} y^{-\beta-\frac 12+\frac 1K-2\varepsilon}
+
M^{-1} y^{-\beta-\frac 12-2\varepsilon}
\biggr).
\end{align*}
But $\frac 12 + \beta < \sigma < 1$, so it is possible to keep
$\varepsilon$ so small that $\beta+1/2-\sigma+2\varepsilon < 0$;
then first term above certainly
tends to $0$ as $M \to \infty$, $y \to 0^+$.
The second and third terms also tend to $0$ so long as we keep
$\varepsilon$ sufficiently small, and $y \geqq M^{-{\minexp}}$;
cf.\ the discussion concerning \eqref{DIOFLEMMA12} above.
\enlargethispage{30pt}
This concludes the proof of Proposition \ref{DIOFLEMMA1}.
\end{Proof}

\begin{proposition} \label{DIOFLEMMA2}
Let $s_0,\beta,\alpha,K,{\minexp}$
be as in Theorem \ref{MAINTHM2}.
Let $v \in \Z$, $k \in \{1,...,\kappa\}$,
and $s \in \frac 12 + \ii[0, \infty)$ or $s=s_0$.
We then have
\begin{equation*} 
\frac 1M \sum_{m=1}^M E_k(m \alpha + {\ii y},s,2v)
\to 0,
\end{equation*}
uniformly as $M \to \infty$, $y \to 0^+$ 
so long as $y \geqq M^{-{\minexp}}$.
\end{proposition}

\begin{Proof}
We take the Fourier expansion of $E_k(z,s,2v)$ to be as in
\eqref{EFOURIEREXP},
and it then follows from the definition of $\beta$ and 
Lemma \ref{EISENST} that
$|c_n| = O(|n|^{\beta+\varepsilon})$ holds for all $n \neq 0$.
The proof of Proposition \ref{DIOFLEMMA2}
is now almost identical to the proof of Proposition \ref{DIOFLEMMA1}.  
\end{Proof}

We are now ready to prove Theorem \ref{MAINTHM2}.
Just as with Theorem \ref{MAINTHM}, we will actually prove a stronger
result, wherein unbounded test functions are allowed:

$\: \:$\textbf{Theorem \ref{MAINTHM2}'.}
\textit{Let $s_0,\beta,\alpha,K,\nu$ be as in Theorem \ref{MAINTHM2}
(on p.\ \pageref{MAINTHM2}),
and let $\gamma \in [0,s_0)$.
Then for any function $f \in B_{\gamma}(\Tqspace)$
and any constant $C_1 > 0$, we have
\begin{equation*}
\frac 1M \sum_{m=1}^{M} f \left( m\alpha + {\ii y} ,0 \right) 
\to \langle f \rangle,
\end{equation*}
uniformly as $M \to \infty$, $y \to 0^+$
so long as $y \geqq C_1 M^{-\minexp}$.}
\vspace{5pt}

{\bf Proof of Theorem \ref{MAINTHM2}'.}
This proof is similar to the arguments given on
pp.\ \pageref{CONDRESSECSTART}--\pageref{MAINPROPPROOF2},
except that the present case is easier.
We will therefore only give an outline of the argument.

It is sufficient to prove \mbox{Theorem \ref{MAINTHM2}'} in the case
$C_1=1$; 
for we may always increase ${\minexp}$ slightly,
keeping \eqref{MAINTHM2ASS} true.

We first prove that there exists a constant $K_2 = 
K_2(\Gamma,s_0,\alpha,{\minexp})>0$ such that
\begin{equation} \label{MAINTHMPR1}
\frac 1M \sum_{m=1}^M \ygg(m \alpha + {\ii y})^{s_0} \leqq K_2,
\qquad \forall M \in \Z^+, \: \: y \in [M^{-{\minexp}},1].
\end{equation}
To this end, we define $G(z)$, $B$ and $G_1(z)$ as in the proof of 
Lemma \ref{CUSPLEMMA}.
It then follows from Proposition \ref{DIOFLEMMA2} (with $s=s_0$ and $v=0$)
that there are numbers $M_0 \in \Z^+$ and $y_0 > 0$ such that
$\frac 1M \sum_{m=1}^M G_1(m \alpha + {\ii y}) < |B|+2$ holds 
whenever $M \geqq M_0$ and $M^{-{\minexp}} \leqq y \leqq y_0$.
Now \eqref{MAINTHMPR1} follows from \eqref{CUSPLEMMA10},
and the fact that $\ygg(z)$ is bounded in the region
$\min(y_0, M_0^{-{\minexp}}) \leqq y \leqq 1$.

Secondly, we claim that Proposition \ref{DIOFLEMMA2} can be sharpened as
follows:
Given $v \in \Z$, $k \in \{1,...,\kappa\}$ and $T>0$,
we have
\begin{equation} \label{MAINTHMPR2}
\frac 1M \sum_{m=1}^M E_k(m \alpha + {\ii y},\sfrac 12 + it ,2v)
\to 0,
\end{equation}
uniformly over all $t \in [0,T]$ as $M \to \infty$, $y \to 0^+$,
so long as $y \geqq M^{-{\minexp}}$.
This is proved by imitating part of the proof of  Proposition \ref{EAPRLEMMA},
using \eqref{MAINTHMPR1} and Proposition \ref{DIOFLEMMA2} 
in the place of Lemma \ref{CUSPLEMMA} and \eqref{MAINPROPASS2}.

The proof is now completed by mimicing the approximation argument on 
pp.\ \pageref{MAINPROPPROOF4}--\pageref{MAINPROPPROOF2},
using \eqref{MAINTHMPR1}, \eqref{MAINTHMPR2} and 
Proposition \ref{DIOFLEMMA1} in the place of
Lemma \ref{CUSPLEMMA},  Proposition \ref{EAPRLEMMA} 
and Proposition \ref{MAASSWVLEMMA}.
$\square$ $\square$ $\square$
\vspace{5pt}

\begin{remark} \label{BESTPOSSIBLE}
The restriction on $\nu$ in \eqref{MAINTHM2ASS} is optimal for
our method of proof in the following precise sense:

\textit{Given any $\beta \geqq 0$, $K \geqq 2$, $v \in \Z$ and
$s \in \frac 12 + i[0,\infty)$ or $s \in (\frac 12,\infty)$,
let $\nu_0$ be the number given in the right hand side of
\eqref{MAINTHM2ASS}.
Then there exists a number $\alpha$ of type $K$, and complex numbers 
$\{ c_n \}_{n \neq 0}$ satisfying $c_n = O(|n|^{\beta})$,
such that the following sum of absolute values
\begin{equation} \label{BESTPOSSIBLE1}
\frac 1M \sum_{n \neq 0} 
\frac {|c_n|}{\sqrt{|n|}} \cdot
\Bigl | W _{v \cdot \sgn(n) ,s - \frac 12} (4 \pi |n| y) \Bigr |
\cdot \biggl | \sum_{m=1}^M e(nm\alpha) \biggr |
\end{equation}
does {\rm not} tend to $0$ as $y = M^{-\nu_0}$, $M \to \infty$.
(For $\beta \leqq \frac 12$ we may even take the $c_n$'s to satisfy
a ``Rankin-Selberg type formula'' 
$\sum_{1 \leqq |n| \leqq N} |c_n|^2 = \text{[const]} \cdot N 
+ O(N^{2\beta})$ as $N \to \infty$.)}

We omit the proof.

Notice that \eqref{BESTPOSSIBLE1} is the sum which we have to treat
if we want to bound the sum \eqref{DIOFLEMMA10} in 
Proposition \ref{DIOFLEMMA1} term by term. 
(The same type of sum arises
also for the Eisenstein series, 
cf.\ Proposition \ref{DIOFLEMMA2}.)
Hence, our remark shows that to obtain any improvement upon 
\eqref{MAINTHM2ASS} in Theorem \ref{MAINTHM2} 
using \eqref{DIOFLEMMA10}, one would have to prove 
cancellation between the terms in \eqref{DIOFLEMMA10}.
\end{remark}

\label{PROOFOFREMARK}
{\bf Proof of the last statement in Remark \ref{SARNAKREMARK}.}
Take $\alpha \in \R$ of type $K \geqq 2$, and 
take $0 < {\minexp} < (K-1)^{-1}$.
We will show how to use Sarnak's theorem on the asymptotic equidistribution
of closed horocycles 
to prove that, for any bounded continuous $\Gamma$-invariant
function $f: \HH \to \C$, and any fixed number $C_1>0$, we have 
\begin{equation} \label{REMPROOF1}
\frac 1M \sum_{m=1}^{M} f \left( m\alpha + {\ii y} ,0 \right) 
\to \langle f \rangle,
\end{equation}
uniformly as $M \to \infty$, $y \to 0^+$
so long as $y \geqq C_1M^{-{\minexp}}$.

By standard approximation arguments, if suffices to prove 
\eqref{REMPROOF1} for $f \in C^{\infty}(T_1\HH) \cap C_c(\Tqspace)$.
For each $y>0$ the function $f(x+{\ii y},0)$ is invariant under 
$x \mapsto x+1$, and hence we have
\begin{align*}
f(x+{\ii y},0) = \sum_{n \in \Z} a(y,n) e(nx),
\end{align*}
where
\begin{align} \label{REMPROOF22}
a(y,n) = \int_0^1 f(x+{\ii y},0) e(-nx) \, dx.
\end{align}
It follows from Sarnak's theorem (cf.\ Theorem \ref{sarnaks} above)
that
\begin{equation} \label{REMPROOF3}
a(y,0) \to \langle f \rangle, \qquad \text{as } \: y \to 0^+.
\end{equation}

Next, we will prove a bound on $a(y,n)$ for $n \neq 0$.
Given a fixed integer $A \geqq 0$, we may apply integration by parts $A$
times in \eqref{REMPROOF22} to obtain
\begin{equation} \label{REMPROOF6}
|a(y,n)| \leqq O ( |n|^{-A}) \sup_{x \in [0,1]} 
\Bigl | \frac{\partial^A}{\partial x^A} f(x+{\ii y},0) \Bigr |.
\end{equation}
But we have, for arbitrary $z=x+iy \in \HH$:
\begin{equation} \label{REMPROOF5}
\Bigl | \frac{\partial^A}{\partial x^A} f(x+{\ii y},0) \Bigr | 
= O(y^{-A})
\end{equation}
(the implied constant depends on $f$, $A$ and $\Gamma$).
To prove \eqref{REMPROOF5}, it is convenient to use the
standard identification of $T_1 \HH$ with the Lie group
$G=\PSL(2,\R)$, given by $G \ni U \mapsto U(i,0) \in T_1 \HH$.
Under this identification, $f$ is a function in
$C^{\infty}(G)$ which is $\Gamma$-left invariant
and which has compact support modulo $\Gamma$.
We let $\minspecX :C^{\infty}(G) \to C^{\infty}(G)$ be the left invariant 
differential operator given by
\begin{equation*}
(\minspecX F)(g) = \frac{d^A}{dt^A} F \bigl( g \matr 1t01 \bigr) \Bigl | _{t=0}.
\end{equation*}
One then verifies that under our identification 
$G \leftrightarrow T_1\HH$, we have
\begin{equation*}
(\minspecX f) \Bigl( \matr {\sqrt{y}}{x/\sqrt{y}}0{1/\sqrt{y}} \Bigr)
= y^A \frac{\partial^A}{\partial x^A} f(x+iy,0).
\end{equation*}
However, since $\minspecX$ is a left invariant operator
and $f$ is a $\Gamma$-left invariant function,
$\minspecX f$ is a  $\Gamma$-left invariant function too. 
Furthermore, $\minspecX f$ has compact support modulo
$\Gamma$, since this is true for $f$. 
Hence $\minspecX f$ is uniformly bounded over all of $G$.
Clearly, \eqref{REMPROOF5} follows from this.

We now obtain, from \eqref{REMPROOF6} and \eqref{REMPROOF5}:
\begin{equation} \label{REMPROOF2}
|a(y,n)| \leqq O ( (|n|y)^{-A}).
\end{equation}

Let us keep $M \in \Z^+$ and $0<y<1$. We have:
\begin{align*}
& \Bigl | \frac 1M \sum_{m=1}^{M} f \left( m\alpha + {\ii y} ,0 \right) 
- a(y,0) \Bigr |
\\
& =
\Bigl | \frac 1M \sum_{n \neq 0} a(y,n) 
\sum_{m=1}^{M} e(m n \alpha) \Bigr |
\leqq
\frac 1M \sum_{n \neq 0} |a(y,n)| \min(M,|| n \alpha ||^{-1}).
\end{align*}
We use the inequality \eqref{REMPROOF2} with $A=0$ for 
$|n| \leqq 1/y$, and with some $A$ yet to be fixed for
$|n| > 1/y$. We then get:
\begin{align*}
\leqq &
O(M^{-1}) \sum_{1 \leqq n \leqq 1/y} \min(M, || n \alpha ||^{-1})
\\
& + 
O(M^{-1}) \sum_{k=1}^{\infty} 
k^{-A}
\sum_{k/y < n \leqq (k+1)/y} \min(M, || n \alpha ||^{-1})
\end{align*}
We now fix $\varepsilon > 0$ so small that $\nu < \frac 1{K-1+\varepsilon}$.
By \eqref{DIOFSUMLEMMARES} in Lemma \ref{DIOFSUMLEMMA}, the above sum is
\begin{align*}
& \leqq O( M^{-1} y^{1-K-\varepsilon})
+ O( M^{-1} ) \sum_{k=1}^{\infty} k^{-A} \Bigl( \frac{k+1}y \Bigr)^{K-1+\varepsilon}.
\end{align*}
If we fix $A$ so large that $K-1+\varepsilon-A < -1$, then we conclude:
\begin{equation} \label{REMPROOF4}
\Bigl | \frac 1M \sum_{m=1}^{M} f \left( m\alpha + {\ii y} ,0 \right) 
- a(y,0) \Bigr |
= O( M^{-1} y^{1-K-\varepsilon}).
\end{equation}
Now \eqref{REMPROOF1} follows from \eqref{REMPROOF3}
and \eqref{REMPROOF4}, since $\minexp < (K-1+\varepsilon)^{-1}$.
$\square$ $\square$ $\square$

\section{Negative results}
\label{MOTEXSECT}

\begin{proposition} \label{NEGCONGRPROP}
Let $\Gamma$ be a subgroup of $\PSL(2,\Z)$ of finite index, and let 
${\minexp}>0$.
Let $\alpha \in \R$ be any irrational number such that there are sequences of 
integers $p_1,p_2,...$ and $0<q_1<q_2<...$ satisfying
\begin{equation*}
\Bigl | \alpha - \frac{p_j}{q_j} \Bigr | \leqq
(3q_j)^{-2-2/{\minexp}}, \qquad \text{for } \: j=1,2,...
\end{equation*}
Then there exist a non-negative function 
$f \in C^{\infty}(T_1 \HH) \cap C_c(\Tqspace)$ 
with $\langle f \rangle = 1$, 
and a sequence of integers $0<M_1<M_2<...$, 
such that
\begin{equation} \label{NEGCONGRRES}
\frac 1{M_j} \sum_{m=1}^{M_j} f(m\alpha + iM_j^{-{\minexp}},0) = 0,
\qquad \text{for all } \: j=1,2,...
\end{equation}
\end{proposition}

\begin{Proof}
Clearly, it suffices to prove the result for $\Gamma = \PSL(2,\Z)$.

Let us consider any fixed $M \geqq 1$, $y>0$ and $j \geqq 1$,
and write $p=p_j,$ $q=q_j$.
We will use the fact that $\alpha \approx p/q$ to show that,
under certain conditions on the sizes of $M$ and $y$,
\textit{all} points $\alpha + iy, 2\alpha+iy,..., M\alpha+iy$
lie far out in the cusp on $\qspace$.

Given $m \in \{1,...,M\}$ we write $d=\gcd(q,mp)$;
then $\gcd(-q/d,mp/d)=1$, 
and hence there is a $\PSL(2,\Z)$-transformation of the form
\begin{equation*}
T = \begin{pmatrix} * & * \\-q/d & mp/d \end{pmatrix}
\in \Gamma = \PSL(2,\Z).
\end{equation*}
We now have
\begin{align*}
\tim T(m\alpha + iy) & = y \Bigl | -\frac qd (m\alpha + iy) + \frac{mp}d \Bigr | ^{-2}
= \frac{d^2y}{q^2y^2 + m^2(p-q\alpha)^2}
\\
& \geqq \frac y {q^2y^2 + M^2q^2(3q)^{-4-4/{\minexp}}}.
\end{align*}
Clearly, this is $\geqq 2$ whenever
\begin{equation} \label{NEGCONGR1}
q^2y^2 \leqq y/4 \quad \text{and} \quad M^2q^2(3q)^{-4-4/{\minexp}} \leqq y/4.
\end{equation}
Hence, for $M,y$ satisfying \eqref{NEGCONGR1},
we have $\ygg(m\alpha+iy) \geqq 2$ for all $m \in \{1,...,M\}$.

Taking $y=M^{-{\minexp}}$, one finds by a quick computation
that \eqref{NEGCONGR1} holds if and only if
\begin{equation} \label{NEGCONGR2}
4^{1/{\minexp}}q^{2/{\minexp}} \leqq M \leqq
4^{-\frac 1{2+{\minexp}}} 
3^{\frac{4(1+{\minexp})}{{\minexp}(2+{\minexp})}} q^{2/{\minexp}}.
\end{equation}
Notice that $4^{1/{\minexp}} 
< 4^{-\frac 1{2+{\minexp}}} 
3^{\frac{4(1+{\minexp})}{{\minexp}(2+{\minexp})}}$,
since $4<3^2$.
Hence, whenever $q$ is sufficiently large, 
there exists at least one integer $M$ 
satisfying \eqref{NEGCONGR2}.
Because of $q=q_j \to \infty$ as $j \to \infty$, we can now
certainly find a sequence of integers
$0<M_1<M_2<M_3<...$ such that for each $j=1,2,3,...$,
and each $m \in \{1,...,M_j\}$,
we have $\ygg(m\alpha+iM_j^{-{\minexp}}) \geqq 2$.

But $R = \{ z \in \HH \mid \ygg(z) < 2 \}$ is an open,
non-empty region on $\qspace$. 
($R$ is open because of the continuity of $\ygg(z)$,
and to see that $R$ is non-empty we need only check that, e.g., 
$\ygg(i)=1$.)
Hence there exists a smooth, $\Gamma$-invariant function 
$f_0 : \HH \to [0,\infty)$
which has compact support contained in $R$, and which is not identically $0$.
We define $f(z,\theta) = f_0(z)$.
Then $f \in C^{\infty}(T_1 \HH) \cap C_c(\Tqspace)$,
and \eqref{NEGCONGRRES} holds.
Rescaling $f$, we can also make $\langle f \rangle = 1$ hold.
\end{Proof}

\section{The pair correlation density of $n^2 \alpha$ mod 1}
\label{APPLSECTION}

The objective of this final section is to establish 
that the equidistribution of the Kronecker sequence $m\alpha$
along closed horocycles implies Poisson statistics for the pair correlation 
density of the sequence $n^2\alpha$ mod 1. 
Rudnick and Sarnak's recent result \cite[Theorem 1 for $d=2$]{RS},
which says that the pair correlation density of $n^2\alpha$ mod 1
is Poissonian for generic $\alpha$ (in Lebesgue measure sense),
is therefore implied by our equidistribution theorem
(Theorem \ref{MAINTHM}).

It would be interesting to see to what extent the convergence
properties of the spacing distributions studied in \cite{RSZ}, 
where $\alpha$ is taken to be
well approximable by rationals, are related to the equidistribution
of Kronecker sequences.

Statistical properties of $n^2\alpha$ mod 1 were also considered 
in connection with the Berry-Tabor conjecture \cite{Berry77} on 
the energy level statistics of integrable Hamiltonian quantum 
systems [the relevant system is here the ``boxed \mbox{oscillator}'' 
with energy levels $n^2\alpha + m$, where $m,n\in\Z_+$], 
integrable quantum maps \cite{Zelditch98}, the ``quantum 
kicked-rotator'' \cite{Casati87,Sinai88,Pellegrinotti88}, and scattering  
problems on certain surfaces of revolution \cite{Zelditch99}.


For any interval $[a,b]$ the {\em pair correlation function} 
is defined as
\begin{equation} \label{PAIRCORR}
R_2([a,b], \alpha, N) = \frac 1N
\Bigl | \Bigl\{
1 \leqq j \neq k \leqq N \: \, \big | \: \, j^2 \alpha - k^2 \alpha 
\in \Bigl [ \frac aN,\frac bN \Bigr ] + \Z  \Bigr\} \Bigr |.
\end{equation}

In the following, we will consider the special
Fuchsian group 
\begin{equation} \label{GAMMA14DEF}
\overline\Gamma_1(4) =\{ (\pm T) \in \PSL(2,\R) \mid T \in \Gamma_1(4) \},
\end{equation}
where $\Gamma_1(4)$ is the congruence subgroup of $\SL(2,\Z)$,
\begin{equation*}
\Gamma_1(4) = \Bigl \{
\begin{pmatrix} a & b \\ c & d \end{pmatrix} \in \SL(2,\Z)
\: \Big | \: c \equiv 0, \, a \equiv d \equiv 1 \pmod 4 \Bigr \}.
\end{equation*}
(Notice that $\overline\Gamma_1(4)$ has a normalized cusp at $\infty$.)


\begin{proposition} \label{APPLPROP}
If $\langle \overline\Gamma_1(4),2,\alpha \rangle$-PSE$_{1/2}$
(cf.\ the definition on p.\ \pageref{PSEDEFPAGE})
holds for some $\alpha \in \R$, then, for any interval $[a,b] \subset \R$:
\begin{equation} \label{APPLPROPRES}
R_2([a,b],\alpha,N) \to b-a, \qquad \text{as } \: N \to \infty .
\end{equation}
That is, the pair correlation function of $n^2\alpha$ mod 1
converges to the one of
independent random variables from a Poisson process.
\end{proposition}

We remark that by Theorem \ref{MAINTHM}',
the assumption in Proposition \ref{APPLPROP} is
satisfied for almost all $\alpha \in \R$, in Lebesgue measure sense.

To prepare for the proof of Proposition \ref{APPLPROP},
let us note that (for arbitrary $\Gamma$ as in earlier sections)
if both
$\langle \Gamma,{\minexp},\alpha \rangle$-PSE$_{\gamma}$
and
$\langle \Gamma,{\minexp},-\alpha \rangle$-PSE$_{\gamma}$
hold, then we also have asymptotic equidistribution of point sets with
an arbitrary weight function $h$, as follows.

\begin{lemma} \label{VIKTFUNKTION}
Let $\Gamma,{\minexp},\alpha,\gamma$ be such that both 
$\langle \Gamma,{\minexp},\alpha \rangle$-PSE$_{\gamma}$
and 
$\langle \Gamma,{\minexp},-\alpha \rangle$-PSE$_{\gamma}$
hold.
Let $h:\R \to \R$ be a bounded, piecewise continuous function of
compact support.
We then have, for any given $f \in B_{\gamma}(\Tqspace)$:
\begin{align} \label{PROPVARIANTRES}
& \frac 1M \sum_{\substack{m \in \Z \\ m \neq 0}} 
h\Bigl( \frac mM \Bigr)
f \left( m\alpha + \ii M^{-{\minexp}} ,0 \right) 
\to  \int_{\R} h(u) \, du \cdot \langle f \rangle,
\end{align}
as $M \to \infty$.
\end{lemma}

\begin{Proof}
It suffices to prove \eqref{PROPVARIANTRES} for functions $h$
which satisfy $h(u)=0$ for all $u < 0$, 
for then functions $h$ 
which satisfy $h(u)=0$ for all $u > 0$
can be treated by replacing
$\alpha$ by $-\alpha$, and the general case follows by
adding the two cases.

If $h(u)= \chi_{(0,b]}$ for some $b>0$, i.e., if $h(u)$ is the
characteristic function of the interval $(0,b]$,
then \eqref{PROPVARIANTRES} follows immediately from
$\langle \Gamma,{\minexp},\alpha \rangle$-PSE$_{\gamma}$
and the fact that $\lim_{M \to \infty} [bM]/M = b$.
Using the relation $\chi_{(b_1,b_2]} = \chi_{(0,b_2]} -
\chi_{(0,b_1]}$,
we now find that \eqref{PROPVARIANTRES} holds for any function
$h$ in the family $\mathcal{R}$ of finite linear
combinations of characteristic functions of intervals
$(b_1,b_2] \subset [0,\infty)$.

Next, if $h(u)$ is a piecewise continuous function on $[0,\infty)$
of compact support, then $h(u)$ is Riemann integrable, and thus there
are sequences of functions $h_i^+,h_i^- \in \mathcal{R}$ such that 
$h_i^-(u) \leqq h(u) \leqq h_i^+(u)$ for all $u > 0$, and 
$\int_0^{\infty} (h_i^+-h_i^-) \, du \to 0$ as $i \to \infty$.
Now if $f \geqq 0$ then
$\sum_m h_i^-(m/M) f(...) \leqq 
\sum_m h(m/M) f(...)
\leqq \sum_m h_i^+(m/M) f(...)$,
and applying \eqref{PROPVARIANTRES} for 
$h_i^+$ and $h_i^-$,
and letting $i \to \infty$, we find that \eqref{PROPVARIANTRES} holds
for $h(u)$.
The case of arbitrary $f$ then follows by 
applying the preceding
result separately to the positive and negative parts of $\text{Re}(f)$ 
and of $\text{Im}(f)$ (each of these four functions belongs to
$B_{\gamma}(\Tqspace)$). 
\end{Proof}

\begin{remark} \label{VREMARK}
Let $T \mapsto \widetilde T$ be the automorphism of $\PSL(2,\R)$
given by $\widetilde T =  \matr i00{-i} T \matr {-i}00i$
(i.e., ${\matr abcd} \: \widetilde{} =  \matr a{-b}{-c}d$).
Concerning the assumption in Lemma \ref{VIKTFUNKTION},
we then have:

\textit{If the lattice $\Gamma$ is invariant under 
$T \mapsto \widetilde T$,
then 
$\langle \Gamma,{\minexp},-\alpha \rangle$-PSE$_{\gamma}$
is equivalent to
$\langle \Gamma,{\minexp},\alpha \rangle$-PSE$_{\gamma}$,
for any given ${\minexp},\alpha,\gamma$.}

In particular, this applies when $\Gamma = \overline\Gamma_1(4)$.

To prove our claim, 
we assume that $\Gamma$ is invariant under $T \mapsto \widetilde T$.
We define $V(z,\theta) := (-\overline{z},-\theta)$;
one then checks that $V \circ T \equiv \widetilde{T} \circ V
: T_1\HH \to T_1\HH$, for all $T \in \PSL(2,\R)$.
Hence, a function $f: T_1\HH \to \C$ is $\Gamma$-invariant 
if and only if $f \circ V$ is $\Gamma$-invariant.
Furthermore, if $\eta$ is a cusp, then so is $-\eta$,
and $\Gamma_{-\eta} = \{ \widetilde{T} \mid T \in \Gamma_{\eta} \}.$
It follows that the set $\minss$ in \eqref{YGGINTRINSIC} is
invariant under $N \mapsto \widetilde{N}$,
and using this one easily checks that
$\ygg(-\overline{z}) = \ygg(z)$ for all $z \in \HH$.

Now our claim follows directly from the definition
of $\langle \Gamma,{\minexp},\alpha \rangle$-PSE$_{\gamma}$
on \mbox{p.\ \pageref{PSEDEFPAGE}},
since replacing $f$ with $f \circ V$ transforms the sum
$\frac 1M \sum_{m=1}^{M} f \left( m\alpha + {\ii y} ,0 \right)$
into
$\frac 1M \sum_{m=1}^{M} f \left( m(-\alpha) + {\ii y} ,0 \right)$.
\end{remark}

{\bf Proof of Proposition \ref{APPLPROP}.}
Fix $\alpha$ as in the proposition.
We will consider the following smoothed version of the
pair correlation function \eqref{PAIRCORR}:
\begin{equation} \label{APPLPROP3}
R_2(g,\psi, \alpha, N) = \frac 1N
\sum_{\substack{j,k \in \Z \\ |j| \neq |k|}} \psi \bigl ( \frac jN \bigr )
\psi \bigl ( \frac kN \bigr )
\sum_{m \in \Z} g(N( j^2 \alpha - k^2 \alpha +m)).
\end{equation}
Here we assume $\psi$ to be an even real-valued
function in $C_c^{\infty} (\R)$, 
and $g$ to be a smooth function such that
$g(x) = O(|x|^{-2})$ as $|x| \to \infty$,
and such that the Fourier transform
$\widehat{g}(u)= \int_{\R} g(x) e(-ux) \, dx$ has compact support.
We will show that for any such functions $\psi$ and $g$, we have
\begin{equation} \label{APPLPROP1}
R_2(g,\psi, \alpha, N) \to \int_{\R} g(x) \, dx \cdot 
\biggl( \int_{\R} \psi(x) \, dx \biggr)^2,
\end{equation}
as $N \to \infty$.

\enlargethispage{10pt}
Let us first prove that \eqref{APPLPROP1} implies
\eqref{APPLPROPRES}. 
Define $\eta(x)$ to be $1$ for $0 < |x| \leqq 1$,
and $0$ for $x=0$ and for $|x| > 1$.
We can now pick two sequences of even non-negative functions
$\psi_i^+$, $\psi_i^- \in C_c^{\infty} (\R)$
with $\psi_i^- \leqq \eta \leqq \psi_i^+$
(hence in particular $\psi_i^-(0)=0$, $\psi_i^+(0) \geqq 1$),
such that
$\int_{\R} (\psi_i^+-\psi_i^-) \, dx \to 0$ as $i \to \infty$. 
Given $[a,b] \subset \R$, we can also pick two sequences of
functions $g_i^+$, $g_i^- \in C^{\infty} (\R)$,
each satisfying $g_i^{\pm}(x) = O_i(|x|^{-2})$ as $|x| \to \infty$
and $\widehat{g_i^{\pm}} \in C_c(\R)$,
such that $g_i^- \leqq \chi_{[a,b]} \leqq g_i^+$
and $\int_{\R} (g_i^+-g_i^-) \, dx \to 0$
as $i \to \infty$.
(Such sequences $g_i^+$, $g_i^-$ can be constructed following
\cite[8.6.1-2]{MarklofMATHANN}.
Cf.\ also \cite{Vaaler}.)
Applying \eqref{APPLPROP1} for $g=g_i^-,$ $\psi=\psi_i^-$
and for $g=g_i^+,$ $\psi=\psi_i^+$,
and then letting $i \to \infty$,
we find that \eqref{APPLPROP1} also holds for
$g=\chi_{[a,b]}$, $\psi=\eta$, i.e., 
\begin{equation} \label{APPLPROP2}
\frac 1N
\sum_{0 < |j| \neq |k| \leqq N}
\sum_{m \in \Z} \chi_{[a,b]}(N(j^2 \alpha - k^2 \alpha +m))
\to 4(b-a),
\end{equation}
as $N \to \infty$.
Notice that if $N > b-a$, then the innermost sum in
\eqref{APPLPROP2} is $1$ if $j^2 \alpha - k^2 \alpha \in [a/N,b/N] + \Z$,
otherwise $0$.
Notice also that $j^2 \alpha - k^2 \alpha$ is invariant under
$j \leftrightarrow -j$ and $k \leftrightarrow -k$.
It follows that for $N > b-a$, the left hand side in \eqref{APPLPROP2} 
equals $4 R_2([a,b],\alpha,N)$.
Hence \eqref{APPLPROP1} implies \eqref{APPLPROPRES}.

It now remains to prove \eqref{APPLPROP1}.
Let $R_2^{(a)}(g,\psi, \alpha, N)$ be the same as in \eqref{APPLPROP3}
but without the restriction $|j| \neq |k|$, i.e., with the
outer sum taken over \textit{all} $j,k \in \Z$.
Since $\psi$ is even, we then have
\begin{align} \label{APPLPROP10}
R_2(g,\psi, \alpha, N) 
=
& R_2^{(a)}(g,\psi, \alpha, N)
- \frac 2N \Bigl( \sum_{j \in \Z} \psi \bigl ( \frac jN \bigr )^2 
-\sfrac 12 \psi(0)^2 \Bigr)
\cdot \sum_{m \in \Z} g(Nm).
\end{align}
Notice here that
\begin{align} \label{APPLPROP11}
\lim_{N \to \infty}
\frac 2N \Bigl( \sum_{j \in \Z} \psi \bigl ( \frac jN \bigr )^2 
-\sfrac 12 \psi(0)^2 \Bigr)
\cdot \sum_{m \in \Z} g(Nm) 
= 2 g(0) \int_{\R} \psi(x)^2 \, dx,
\end{align}
since $\psi \in C_c^{\infty}(\R)$ and $g(x) = O(|x|^{-2})$.
Also, by the Poisson summation formula, we have
\begin{align*}
R_2^{(a)}(g,\psi, \alpha, N) 
& = \frac 1{N^2} \sum_{j,k \in \Z} \psi \bigl ( \frac jN \bigr )
\psi \bigl ( \frac kN \bigr )
\sum_{m \in \Z} \widehat{g} \Bigl( \frac mN \Bigr)
e(m (j^2-k^2) \alpha )
\\
& = \frac 1{N} \sum_{m \in \Z} \widehat{g} \Bigl( \frac mN \Bigr) 
\Bigl | \theta_{\psi}(m \alpha + \ii N^{-2} ) \Bigr |^2,
\end{align*}
where we have defined
\begin{equation*}
\theta_{\psi}(x+{\ii y}) = 
y^{1/4} \sum_{j \in \Z} \psi \bigl ( jy^{1/2} \bigr ) e(j^2x).
\end{equation*}
This\enlargethispage{10pt} is a theta sum with a smooth cutoff function;
such sums were studied in \cite{Marklof}
(both for smooth and sharp cutoff functions).
To apply \cite{Marklof} to our situation,
we take $\Gamma = \overline{\Gamma}_1(4)$ (cf.\ \eqref{GAMMA14DEF}) and
$f=|\Theta_{\psi}|^2$, where $\Theta_{\psi}$ is as in
\cite[Prop.\ 3.1]{Marklof}.
Notice that $\Theta_{\psi}$ is a function on a 
\textit{$4$-fold cover} of $\PSLqspace$. 
However, we claim that
$|\Theta_{\psi}|^2$ is a function on $\PSLqspace$ itself.
To see this, we first use \cite[(10),(12)]{Marklof} 
to check that the function
$\theta_k$ in \cite[(24)]{Marklof} satisfies
$\theta_k \bigl( [\matr {-1}00{-1}, \beta_{-1}] (z,\phi) \bigr)
= e^{-i(k+\frac 12) \pi} \theta_k(z,\phi)$.
We also notice that for $k \in \Z^+$ odd, $\theta_k \equiv 0$.
Hence $\Theta_{\psi}([\matr {-1}00{-1}, \beta_{-1}] (z,\phi))
= e^{-i \pi/2} \Theta_{\psi}(z,\phi)$, by \cite[(26)]{Marklof},
and thus $|\Theta_{\psi}(z,\phi)|^2$ is invariant under
$[\matr {-1}00{-1}, \beta_{-1}]$.
This gives the desired result, since
we already know that $|\Theta_{\psi}(z,\phi)|^2$ is invariant under 
$\Delta_1(4)$, cf.\ \cite[(15)]{Marklof}.
In conclusion, we see from
\cite[Prop.\ 3.1]{Marklof}
that the function $f=|\Theta_{\psi}|^2$
is smooth and $\Gamma$-invariant, and 
\begin{equation*}
f(z,0) = |\theta_{\psi}(z)|^2, 
\qquad \forall z \in \HH.
\end{equation*}

It follows from \cite[Prop.\ 3.2]{Marklof} that $f \in
B_{\frac 12}(\Tqspace)$.
Hence, since $\langle \overline\Gamma_1(4),2,\alpha \rangle$-PSE$_{1/2}$
holds by assumption, Lemma \ref{VIKTFUNKTION}
and Remark \ref{VREMARK} apply, and we obtain
\begin{align} \label{APPLPROP31}
\lim_{N \to \infty} \Bigl( R_2^{(a)}(g,\psi, \alpha, N) 
- \frac 1N \widehat{g} (0) \bigl | \theta_{\psi}( \ii N^{-2} ) 
\bigr |^2 \Bigr) 
& =
\int_{\R} \widehat{g} (u) \, du \cdot \langle f \rangle.
\end{align}
But here $\int_{\R} \widehat{g} (u) \, du = g(0)$,
and it follows from \cite[(48)]{Marklof}
and $\mu(\qspace) = 2\pi$ 
that $\langle f \rangle = 2 \int_{\R} \psi(x)^2 \, dx$.
Finally, the definition of $\theta_{\psi}(x+{\ii y})$ implies that
\begin{align} \label{APPLPROP30}
\frac 1{\sqrt N} \theta_{\psi}( \ii N^{-2} )
= \frac 1N \sum_{j \in \Z} \psi(j/N) \to \int_{\R} \psi(x) \, dx,
\qquad
\text{as } \: N \to \infty.
\end{align}
Now \eqref{APPLPROP1} follows from 
\eqref{APPLPROP10}, \eqref{APPLPROP11},
\eqref{APPLPROP31} and
\eqref{APPLPROP30}, and the proof is complete.
$\square$ 
\vspace{5pt}

\setlength{\parskip}{0pt}

\vspace{5pt}

\begin{scriptsize}
\noindent
\textsc{Jens Marklof,} 
School of Mathematics, University of Bristol,
Bristol BS8 1TW, U.K.

\vspace{5pt}

\noindent
\textsc{Andreas Str\"ombergsson,}
Department of Mathematics, Princeton University,
Fine Hall,\linebreak Washington Road, Princeton, NJ 08544, U.S.A.
\end{scriptsize}


\begin{thebibliography}{MMMM}

\addtocounter{bibno}{1}
\bibitem[AS]{AS}
M. Abramowitz and I. Stegun,
\textit{Handbook of Mathematical Functions},
AMS \textbf{55}, 7th ed., 1968.

\addtocounter{bibno}{1}
\bibitem[BR]{BR}
J. Bernstein and A Reznikov,
Analytic continuation of representations and estimates of automorphic
forms,
{\em Ann.\ of Math.}, \textbf{150} (1999), 329--352.

\addtocounter{bibno}{1}
\bibitem[BT]{Berry77}
M.V.\,Berry and M.\,Tabor, 
Level clustering in the regular spectrum,
{\em Proc. Roy. Soc.} A {\bf 356} (1977) 375-394.

\addtocounter{bibno}{1}
\bibitem[CGI]{Casati87}
G. Casati, I. Guarneri and F. M. Izrailev,
Statistical Properties of the Quasi-Energy Spectrum of a Simple
Integrable System, {\em Phys. Lett.} A \textbf{124} (1987), 263--266. 

\addtocounter{bibno}{1}
\bibitem[D1]{Deligne}
P. Deligne, La conjecture de Weil I, II, {\em Publ. IHES} 
\textbf{43} (1974), 273--307;
\textbf{52} (1980), 137--252.

\addtocounter{bibno}{1}
\bibitem[D2]{Dunster1}
T. M. Dunster, Uniform asymptotic expansions for Whittaker's
confluent hypergeometric functions, {\em SIAM J. Math. Anal.} \textbf{20}
(1989), 744--760.

\addtocounter{bibno}{1}
\bibitem[D3]{Dunster2}
T. M. Dunster, Uniform asymptotic approximations for the
Whittaker functions $M_{\kappa,i\mu}(z)$
and $W_{\kappa,i\mu}(z)$,
to appear in {\em Analysis and Applications}.

\addtocounter{bibno}{1}
\bibitem[EMM]{Eskin98}
A.\,Eskin, G.\,Margulis and S.\,Mozes, Upper bounds and
asymptotics in a quantitative version of the Oppenheim conjecture,
{\em Ann. of Math.} {\bf 147} (1998) 93-141.

\addtocounter{bibno}{1}
\bibitem[G]{Good}
A. Good,
Cusp forms and eigenfunctions of the Laplacian,
{\em Math.\ Ann.}\ \textbf{255} (1981), pp.\ 523 -- 548.

\bibitem[H1]{TR1}
D. A. Hejhal,
\textit{The Selberg Trace Formula for PSL$(2,{\R})$}, Vol. 1,
Lecture Notes in Math. \textbf{548}, 
Springer-Verlag, Berlin, 1976.

\addtocounter{bibno}{1}
\bibitem[H2]{TR2}
D. A. Hejhal,
\textit{The Selberg Trace Formula for PSL$(2,{\R})$}, Vol.\ 2,
Lecture Notes in Math. \textbf{1001}, 
Springer-Verlag, Berlin, 1983.

\addtocounter{bibno}{1}
\bibitem[H3]{DHHORO}
D. A. Hejhal,
On value distribution properties of automorphic functions along closed
horocycles,
in \textit{XVI-th Rolf Nevanlinna Colloquium,}
(edited by I. Laine and O. Martio), deGruyter, 1996, pp.\ 39--52.

\addtocounter{bibno}{1}
\bibitem[I1]{IW}
H. Iwaniec, 
\textit{Introduction to the Spectral Theory of Automorphic Forms,}
Biblioteca de la Revista Matem\'atica Iberoamericana,
Madrid, 1995.

\addtocounter{bibno}{1}
\bibitem[I2]{IW2}
H. Iwaniec, 
\textit{Topics in Classical Automorphic Forms,}
Graduate Studies in Mathematics \textbf{17},
AMS, 1997.

\addtocounter{bibno}{1}
\bibitem[J1]{DJ94}
D. Jakobson,
Quantum unique ergodicity for Eisenstein series on
PSL$_2(\Z) \backslash$PSL$_2(\R)$,
{\em Ann. Inst. Fourier} (Grenoble), \textbf{44} (1994), 1477--1504.

\addtocounter{bibno}{1}
\bibitem[J2]{DJ97}
D. Jakobson,
Equidistribution of cusp forms on PSL$_2(\Z) \backslash$PSL$_2(\R)$,
{\em Ann. Inst. Fourier} (Grenoble), \textbf{47} (1997), 967--984.


\addtocounter{bibno}{1}
\bibitem[KS]{KS}
H. H. Kim and P. Sarnak,
Refined estimates towards the Ramanujan and Selberg conjectures,
appendix to 
H. H. Kim, Functoriality for the exterior square of $GL_4$ and symmetric
fourth of $GL_2$, preprint.


\addtocounter{bibno}{1}
\bibitem[M1]{Marklof}
J. Marklof,
Limit theorems for theta sums,
\textit{Duke Math.\ J.\ } \textbf{97} (1999), 127--153.

\addtocounter{bibno}{1}
\bibitem[M2]{MarklofMATHANN}
J. Marklof,
Pair correlation densities of inhomogeneous quadratic forms,
preprint 2000/2002, to appear in {\em Ann.\ of Math.}

\addtocounter{bibno}{1}
\bibitem[N]{Nathanson}
M. B. Nathanson,
\textit{Additive Number Theory},
Springer-Verlag, 1996.

\addtocounter{bibno}{1}
\bibitem[O1]{Olver}
F. W. J. Olver,
\textit{Asymptotics and special functions},
Academic Press, New York and London, 1974.

\addtocounter{bibno}{1}
\bibitem[O2]{OlverW}
F. W. J. Olver,
Whittaker functions with both parameters large: uniform approximation
in terms of parabolic cylinder functions,
{\em Proc. Roy. Soc. Edinburgh Sect.} A \textbf{86} (1980), 213--234.

\addtocounter{bibno}{1}
\bibitem[P]{Pellegrinotti88}
A. Pellegrinotti, 
Evidence for the Poisson distribution for quasi-energies 
in the quantum kicked-rotator model,
{\em  J. Statist. Phys.} {\bf 53} (1988) 1327--1336.

\addtocounter{bibno}{1}
\bibitem[R1]{Ratner91}
M.\,Ratner, On Raghunathan's measure conjecture,
{\em Ann. of Math.} {\bf 134} (1991) 545-607.

\addtocounter{bibno}{1}
\bibitem[R2]{Ratner91b}
M.\,Ratner, Raghunathan's topological conjecture 
and distributions of unipotent flows, {\em Duke Math. J.} {\bf  63}
(1991) 235-280.

\addtocounter{bibno}{1}
\bibitem[RS]{RS}
Z. Rudnick and P. Sarnak,
The pair correlation function of fractional parts of polynomials,
{\em Comm. Math. Phys.} \textbf{194} (1998), 61--70.

\addtocounter{bibno}{1}
\bibitem[RSZ]{RSZ}
Z. Rudnick, P. Sarnak and A. Zaharescu,
The distribution of spacings between the fractional parts of $n^2 \alpha$,
{\em Invent. Math.} \textbf{145} (2001), 37--57.

\addtocounter{bibno}{1}
\bibitem[S1]{Sarnak}
P. Sarnak,
Asymptotic behavior of periodic orbits of the horocycle flow and
Eisenstein series,
{\em Comm. Pure Appl. Math.} \textbf{34} (1981), 719--739.

\addtocounter{bibno}{1}
\bibitem[S2]{Schmidt}
W. M. Schmidt,
\textit{Approximation to algebraic numbers},
{\em S\'erie des Conf\'erences de l'Union Math\'ematique Internationale},
No.\ 2. Monographie No.\ 19 de l'Enseignement Math\'ematique,
Geneva, 1972.
(Also in {\em Enseignement Math.}\ (2) \textbf{17} (1971), 187--253.)



\addtocounter{bibno}{1}
\bibitem[S4]{Sinai88}
Ya. G. Sinai, 
The absence of the Poisson distribution for spacings 
between quasi-energies in the quantum kicked-rotator model. 
{\em Phys.} D {\bf 33} (1988) 314--316.



\addtocounter{bibno}{1}
\bibitem[S6]{Strhoro}
A. Str\"ombergsson,
Some results on the uniform equidistribution of long closed
horocycles,
in \textit{Studies in the analytic and spectral theory of automorphic
forms}, Thesis, Uppsala University, 2001, 137--226.
[Available electronically at: 
http://www.math.princeton.edu/${}^\sim$astrombe/papers.html]

\addtocounter{bibno}{1}
\bibitem[V]{Vaaler}
J. D. Vaaler, Some extremal functions in Fourier analysis,
{\em Bull. Amer. Math. Soc.} {\bf 12} (1985) 183--216.

\addtocounter{bibno}{1}
\bibitem[Z]{Zelditch98}
S. Zelditch, Level spacings for integrable quantum maps in genus zero,
{\em Comm. Math. Phys.} {\bf 196} (1998) 289--318,
Addendum: "Level spacings for integrable quantum maps in genus zero",
{\em ibid.}, 319--329.

\addtocounter{bibno}{1}
\bibitem[ZZ]{Zelditch99}
S. Zelditch and M. Zworski, 
Spacing between phase shifts in a simple scattering problem,
{\em Comm. Math. Phys.} {\bf 204} (1999) 709--729.
\end{thebibliography}
\end{document}